\documentclass[USenglish]{article}	

\usepackage[utf8]{inputenc}				
\usepackage[big,online]{dgruyter}	
\usepackage{lmodern}
\usepackage{microtype}
\usepackage[numbers,square,sort&compress]{natbib}
\usepackage{lineno,hyperref}
\usepackage{caption}
\usepackage{mathrsfs}
\usepackage{makecell}
\usepackage{multirow}
\usepackage{graphics}
\usepackage{subfigure}
\usepackage{url}
\usepackage{bm}
\usepackage{dsfont}
\usepackage[normalem]{ulem}
\usepackage{epstopdf}
\usepackage{graphicx,color}
\usepackage{setspace,epsfig}
\usepackage{algorithm}
\usepackage{supertabular}
\usepackage{longtable}
\usepackage{colortbl}
\usepackage{arydshln}
\usepackage{multirow}
\usepackage{multicol}
\usepackage{algorithmic}
\usepackage{caption}

\theoremstyle{dgthm}
\newtheorem{theorem}{Theorem}

\newtheorem{proposition}{Proposition}

\theoremstyle{dgdef}
\newtheorem{definition}{Definition}

\begin{document}

	\articletype{Research Article}
	\received{Month	DD, YYYY}
	\revised{Month	DD, YYYY}
  \accepted{Month	DD, YYYY}
  \journalname{De~Gruyter~Journal}
  \journalyear{YYYY}
  \journalvolume{XX}
  \journalissue{X}
  \startpage{1}
  \aop
  \DOI{10.1515/sample-YYYY-XXXX}

\title{An adaptive variational model for multireference alignment with mixed noise}
\runningtitle{Short title}

\author[1]{Cuicui Zhao}
\author*[2]{Jun Liu}
\author*[1]{Xinqi Gong}
\affil[1]{\protect\raggedright
Mathematics Intelligence Application LAB, Institute for Mathematical Sciences, Renmin University of China, Beijing 100872, China, e-mail: cuicuizhao@ruc.edu.cn, xinqigong@ruc.edu.cn.}
\affil[2]{\protect\raggedright
Laboratory of Mathematics and Complex Systems (Ministry of Education of China), School of Mathematical Sciences, Beijing Normal University, Beijing, 100875, China, email: jliu@bnu.edu.cn.}
	
	
\abstract{ Multireference alignment (MRA) problem is to estimate an underlying signal from a large number of noisy circularly-shifted observations. The existing methods are always proposed under the hypothesis of a single Gaussian noise. However, the hypothesis of a single-type noise is inefficient for solving practical problems like single particle cryo-EM. In this paper, We focus on the MRA problem under the assumption of Gaussian mixture noise. We derive an adaptive variational model by combining maximum a posteriori (MAP) estimation and soft-max method. There are two adaptive weights which are for detecting cyclical shifts and types of noise. Furthermore, we provide a statistical interpretation of our model by using expectation-maximization(EM) algorithm. The existence of a minimizer is mathematically proved. The numerical results show that the proposed model has a more impressive performance than the existing methods when one Gaussian noise is large and the other is small.}


\maketitle

\section{Introduction}\label{intro}
Multireference alignment(MRA) is the problem of estimating a true signal from a number of noisy and circularly-shifted observations. This mathematical model arises in many scientific and engineering fields, for instance, structural biology\cite{1992Onfdsjpo,2014sdhoisAn,201dshois1A,1992cxhoisShift,2005dssdhoiiMaximum,2012sdhoisOptimal},
single cell genomic sequencing\cite{2015wewNanoOK}, radar\cite{2005shdospUsing,2003cshoiuoiFast},
robotics\cite{rosen2017guiysesync}, crystalline simulations\cite{2013sdhiosoNoisy},image registration and super-resolution\cite{1992dsuhsA,2002Extension,2008oiuoiOptimal} and some algorithms and theoretical analysis\cite{2017Bispectrum,2019dsguThe,2014sdsMultireference,2017csdssds,8352518,
8590822,2020Multi,2021Aizenbudfshois,Romanov_2021,2020dsdsdsWavelet}. There are some variants of MRA problem, such as
heterogeneous MRA\cite{2018Hetefdsohds}, super-resolution MRA\cite{a2021ndorye}.
\par
The mathematical description of MRA problem is
\begin{equation}\label{a1.1}
\pmb{f}_i=R_{l_i}\pmb{u}+\pmb{v}_i,~~~i=1,2,...,M,
\end{equation}
where $\pmb{u}=(u_1,u_2,...,u_N),~\pmb{f}_i=(f_{i,1},f_{i,2},...,f_{i,N}),~
\pmb{v}_i=(v_{i,1},v_{i,2},...,v_{i,N})$ are the real signal, the $i$-th observation
and the corresponding noise respectively . $R_l$ is a circularly shifted operator,
namely, $R_l\pmb{u}[j]=\pmb{u}[j-l]$ in the sense of zero-based indexing and modulo $N$,
where $l\in\{0,1,2,...,N-1\}$. Significantly, both the true signal $x$ and the shifts $r_i$ are unknown.
The goal is to recover $x$ from
these observations $\pmb{f}_i$.
\par
As far as we know, the existing literatures for MRA problem have been always based on the hypothesis
of a single Gaussian noise. The literature can roughly be divided into two patterns. One is
first estimating shifts and then estimating the true signal\cite{2014sdsMultireference}. the other is estimating the true signal
without seeking for shifts\cite{2017Bispectrum,8352518}.
Here we are only concentrated on the latter. For the latter pattern,
there are usually two general approaches. One is based on statistical knowledge
like maximum likelihood estimation(MLE), maximum a posteriori(MAP) estimation, expectation-maximization(EM) algorithm\cite{katsevich2020likelihood,2012RELION}. The other is based on shift variant features such as spectral method largest spectral gap\cite{8352518}, frequency marching(FM)\cite{2017Bispectrum}, optimization on phase manifold\cite{2017Bispectrum} and optimization on phase synchronization\cite{2017Bispectrum}. The advantage of the former has higher accuracy, but the latter
needs less time and computer resources.
\par
MRA problem is a simplified mathematical model for single
particle cryo-electron microscopy (cryo-EM)\cite{2020Single,2019MATHEMATICS,2012RELION}, which is a very
popular technique for visualizing biological molecules \cite{2000Single,2021articldssgiuaye,202daaasd0Single}. Image denoising and image alignment are two tasks in single particle cryo-EM. 
The MRA problem is a simplified model for tackling these two tasks simultaneously. However, in practice, the noise in cryo-EM images\cite{0fsoodisTopaz} are always more complicated beyond a single Gaussian noise.
 Motivated by the above discussions, we extend the type of noise from a single Gaussian noise to Gaussian mixture noise to close the gap between the MRA model and the actual applications.
\par
In this work, we derive an adaptive variational model for MRA problem with Gaussian mixture noise
by combing MAP estimation and soft-max method. The big challenge is that both circularly-shifted translations and the types of
noise are unknown. To solve this difficulty, the proposed model contains two weights, where one is
for the circularly-shifted translation of each observation, the other is for the type of noise on each component of each observation. In addition, we provide a statistical
explanation for our proposed model by using exception-maximization(EM) algorithm. Furthermore,
we prove the existence of a minimizer in our proposed model with total-variation(TV) regularizer.
We design an algorithm to by the alternating direction iterative method and the augmented Lagrange method. In addition, we provide some convergence analysis. In numerical experiments, our proposed model outperforms the existing algorithms when one Gaussian noise is large and the other is small.
\par
{\bf Organization of this paper.} We derive the proposed model in Section \ref{Sec2}, including model hypothesis and modeling process. In Section \ref{Sec3}, we provide a statistical interpretation of the proposed model. We prove the existence of a minimizer in Section \ref{Sec4}. In Section \ref{Sec5}, we design an algorithm by alternating direction iterative method and augmented Lagrange method. Section \ref{Sec6} provides some convergence analysis. Section \ref{Sec7} shows some numerical results. We summarize this paper in Section \ref{Sec8}.
\par
One can reproduce this work by the code in \url{https://github.com/MIALAB-RUC/MRA-MGG-softmax}.
\section{The proposed model}\label{Sec2}
In this section, we derive an adaptive variational model for multireference alignment(MRA) problem with mixed Gaussian-Gaussian(MGG) noise by combining maximum a posteriori(MAP) estimation and soft-max method.
\par
We denote $\mathscr{U},~\mathscr{F},~\mathscr{V},~\mathscr{L}$ as random variables of pixels in the true signal, observation, noise, and circularly-shifted operator respectively. $f,~u,~v,~l$ are corresponding sample values. $P_{\mathscr{X}}(x)$ and $p_{\mathscr{X}}(x)$ represent the cumulative distribution function and the probability density distribution function of random variable $\mathscr{X}$ at point $x$ separately.
\subsection{Model hypothesis}
There are some basic hypotheses in the proposed model:\\
\textbf{A1}~~The true signal is corrupted by some mixed noise with mean $0$. We denote $\pmb{\alpha}=(\alpha_1,\alpha_2,...,\alpha_K)$ as the mixed ratio, where $\alpha_k~$ is the ratio of the $k-$th kind of noise among the mixed noise. Moreover, $\pmb{\sigma^2}=(\sigma_1^2,~\sigma_2^2,...,\sigma_K^2)$ denotes the mixed noise parameters, where $\sigma_k^2$ is the variance of the $k-$th kind of noise. Denote $\pmb{\theta}=(\pmb{\alpha},\pmb{\sigma^2})$ as all parameters of the mixed noise ; \\
\textbf{A2}~~ The value of noise is a realization of random variable $\mathscr{V}$. Each components of observations are mutually independently and identically distributed(i.i.d) with probability density function $p_{\mathscr{V}}$; \\
\textbf{A3}~~The value of a circularly shifted operator is a realization of random variable $\mathscr{L}$. The circularly shifted transformations of different observations are mutually independently and identically distributed(i.i.d) with discrete uniform distribution in $\{0,1,...,N-1\}$; \\
\textbf{A4}~~The true signal follows a Gibbs prior distribution.
\par
By \textbf{A1}, The mixed noise can be expressed as
\begin{equation}\label{a1.2}
  \mathscr{V}=\left\{
   \begin{aligned}
   \mathscr{V}_1,& ~~~\text{when event~$A_1$~occurs},\\
      \mathscr{V}_2, &~~~\text{when event~$A_2$~occurs},\\
      ...&\\
      \mathscr{V}_K,& ~~~\text{when event~$A_K$~occurs},\\
   \end{aligned}
   \right.
  \end{equation}
  where $P(A_k)=\alpha_k$, and $\displaystyle\sum_{k=1}^K\alpha_k=1$. Moreover, we denote $p_{\mathscr{V}_k}(v)$ as the  probability density function of random variable $\mathscr{V}_k$ for any $k=1,2$.
\par
By \textbf{A3}, we can get
\begin{equation}\label{a1.2_1}
  \mathscr{L}=\left\{
   \begin{aligned}
   0,& ~~~\text{when event~$B_0$~occurs},\\
     1, &~~~\text{when event~$B_1$~occurs},\\
      ...&\\
      N-1,& ~~~\text{when event~$B_{N-1}$~occurs},,\\
   \end{aligned}
   \right.
  \end{equation}
where $P(B_l)=\frac{1}{N}$ for any $l=0,1,...,N-1$. Furthermore, we can get
\begin{equation}\label{a1.5.1.1}
  P_{\mathscr{L}}(l)=\frac{1}{N}.
\end{equation}
\par
By \textbf{A4}, the probability density function of the true signal is
\begin{equation}\label{a1.6}
    p_{\mathscr{U}}(u)=\frac{1}{T}e^{-\gamma \phi(u)},
\end{equation}
where $T>0$ is a constant and~$\phi$~is a given function.
\subsection{ The modeling process}\label{subsec_modeling}

In this subsection, we apply MAP estimation and soft-max method to derive an adaptive variational model for MRA problem under the assumption of mixture Gaussian noise.
\\
\subsubsection{MAP estimation}
The observations $\pmb{f}_1,~\pmb{f}_2,...,~\pmb{f}_M$ are some realizations of random variables $\pmb{\mathscr{F}}_1,\pmb{\mathscr{F}}_2,...,\pmb{\mathscr{F}}_M$ respectively. The $i$-th observation is $\pmb{f}_i=(f_{i,1},f_{i,2},...,f_{i,N})$ and its corresponding random variable is $\pmb{\mathscr{F}}_i=(\mathscr{F}_{i,1},\mathscr{F}_{i,2},...,\mathscr{F}_{i,N})$ for any i=1,2,...,M. Furthermore, we denote $\mathscr{A}=(\pmb{\mathscr{F}}_1,\pmb{\mathscr{F}}_2,...,\pmb{\mathscr{F}}_M)$. $A=(\pmb{f}_1,~\pmb{f}_2,...,~\pmb{f}_M)$ is denoted as a realization of $\mathscr{A}$. The true signal $\pmb{u}=(u_1,u_2,...,u_N)$ is a realization of random vector $\pmb{\mathscr{U}}=(\mathscr{U}_1,~\mathscr{U}_2,...,\mathscr{U}_N)$.
The goal is to estimate the true signal $\pmb{u}$ from all observations. We need to maximize $P_{\pmb{\mathscr{U}}|\mathscr{A}}(\pmb{u}|A)$.
\par
\begin{proposition}\label{pro1}
$\mathscr{V},~\mathscr{V}_1,~\mathscr{V}_2,...,~\mathscr{V}_K$~are random variables that satisfy \eqref{a1.2}.
$p_{\mathscr{V}_k}(v)$ is the probability density function of $\mathscr{V}_k$ for any $k=1,2,...,K$; then
\begin{equation}\label{a2.2.4}
p_{\mathscr{V}}(v)=\sum_{k=1}^K\alpha_k p_{\mathscr{V}_k}(v_k).
\end{equation}
\end{proposition}
\begin{proposition}\label{pro2}
Assume that $\pmb{\mathscr{F}}=R_{l}\pmb{\mathscr{U}}+\pmb{\mathscr{V}}$, where $l$ is a known constant. $\pmb{\mathscr{U}}$ and $\pmb{\mathscr{V}}$ are mutually independent. The probability density function of $\pmb{\mathscr{V}}$ is $p_{\pmb{\mathscr{V}}}(\pmb{v})$; then
\begin{equation}\label{a2.2.5}
p_{\pmb{\mathscr{F}}|R_{l}\pmb{\mathscr{U}}}(\pmb f|R_{l}\pmb u)=p_{\pmb{\mathscr{V}}}(\pmb f-R_{l}\pmb u).
\end{equation}
\end{proposition}
\par
By Bayes' law, we can get
\begin{equation}\label{a2.2.6}
  P_{\pmb{\mathscr{U}}|\mathscr{A}}(\pmb{u}|A)=
  \frac{P_{\mathscr{A}|\pmb{\mathscr{U}}}(A|\pmb{u})P_{\pmb{\mathscr{U}}}(\pmb{u})}
  {P_{\mathscr{A}}(A)}.
\end{equation}
 Note that $P_{\mathscr{A}}(A)$ is a constant. By taking the logarithm of \eqref{a2.2.6}, the goal can be converted to the following problem
\begin{equation}\label{a2.2.7}
\max_{\pmb{u},\pmb{\theta}}\left \{\log P_{\mathscr{A}|\pmb{\mathscr{U}}}(A|\pmb{u})+\log P_{\pmb{\mathscr{U}}}(\pmb{u})\right\},
\end{equation}
where $\pmb{\theta}$ denotes all parameters of the mixed noise. Because the random vectors $\pmb{\mathscr{F}}_1,\pmb{\mathscr{F}}_2,...,\pmb{\mathscr{F}}_M$ are mutually independent, we can get
\begin{equation}\label{a2.2.8}
P_{\mathscr{A}|\pmb{\mathscr{U}}}(A|\pmb{u})=\prod_{i=1}^MP_{\pmb{\mathscr{F}}_i|\pmb{\mathscr{U}}}(
\pmb{f}_i|\pmb{u}).
\end{equation}
Substitute \eqref{a2.2.8} into \eqref{a2.2.7}, then the problem \eqref{a2.2.7} equals to the following problem
\begin{equation}\label{a2.2.9}
 \min_{\pmb{u},\pmb{\theta}}\left\{-\sum_{i=1}^M\log P_{\pmb{\mathscr{F}}_i|\pmb{\mathscr{U}}}(
\pmb{f}_i|\pmb{u})-\log P_{\pmb{\mathscr{U}}}(\pmb{u})\right\}.
\end{equation}
By the law of total probability and substituting \eqref{a1.5.1.1} to \eqref{a2.2.9} , we can get
\begin{equation}\label{a2.2.10}
P_{\pmb{\mathscr{F}}_i|\pmb{\mathscr{U}}}(\pmb{f}_i|\pmb{u})=
\sum_{l=1}^NP_{\mathscr{L}}(l)P_{\pmb{\mathscr{F}}_i|R_l\pmb{\mathscr{U}}}(\pmb{f}_i|R_l\pmb{u})
=
\frac{1}{N}\sum_{l=1}^NP_{\pmb{\mathscr{F}}_i|R_l\pmb{\mathscr{U}}}(\pmb{f}_i|R_l\pmb{u}),
\end{equation}
for any $i=1,2,...,M$. By Proposition \ref{pro2}, we can get
\begin{equation}\label{a2.2.11}
P_{\pmb{\mathscr{F}}_i|\pmb{\mathscr{U}}}(\pmb{f}_i|\pmb{u})=
\frac{1}{N}\sum_{l=1}^NP_{\pmb{ \mathscr{V}}}(\pmb{f}_i-R_l\pmb{u}).
\end{equation}
Substitute \eqref{a2.2.11} into \eqref{a2.2.9}, then the problem can be converted to
\begin{equation}\label{a2.2.12}
\min_{\pmb{u},\pmb{\theta}}\left\{\mathcal{L}(\pmb{u},\pmb{\theta})=-\sum_{i=1}^M\log \sum_{l=1}^NP_{\pmb{ \mathscr{V}}}(\pmb{f}_i-R_l\pmb{u})
 -\log P_{\pmb{\mathscr{U}}}(\pmb{u})\right\}.
\end{equation}
We note that there is a logarithm of summation term, which is difficult to handle in energy minimization problem. Next, we will apply soft-max method to solve this difficulty.
\\

\subsubsection{soft-max method}\label{soft-max-section}
 \begin{definition}[Soft-max]\cite{Polyak2016}\label{def1}
Given a vector $\mathbf{x}=(x_{1},x_{2},..x_{N})$; then for any fixed $\varepsilon>0$, the soft-max operator is defined by
\begin{equation}\label{a2.2.13}
\mbox{max}_{\varepsilon}(\mathbf{x}):=\varepsilon\log\sum\limits_{l=1}^{N}\displaystyle\biggl. e^{\frac{x_{l}}{\varepsilon}}.
\end{equation}
\end{definition}
It is easy to check that $\underset{\varepsilon\rightarrow 0}{\lim}~\max_{\varepsilon}(\textbf{x})=\max\{\mathbf{x}\}$.
\begin{proposition}\label{pro3}\cite{Polyak2016}
Set $G_{\varepsilon}(\mathbf{x})=\max_{\varepsilon}(\mathbf{x})$;~then for any fixed $\varepsilon>0$,  $G_{\varepsilon}(\mathbf{x})$ is convex with respect to $\mathbf{x}$.
\end{proposition}
\begin{definition}[Fenchel-Legendre transformation]\label{def2}
\cite{Polyak2016}
~Denote $G^{*}$ the Fenchel-Legendre transformation of $G$, which is defined by
\begin{equation}\label{a2.2.14}
G^*(\mathbf{w}):=\max\limits_{\mathbf{x}}\{<\mathbf{x},\mathbf{w}>-
G(\mathbf{x})\}.
\end{equation}
\end{definition}
\par
\begin{proposition}\label{pro4}\cite{2005dsuhsVariational}
A function $G:\mathbb{R}^M\rightarrow \mathbb{R}\cup\{+\infty\}$ is convex and lower semi-continuity if and only if $G=G^{**}$.
\end{proposition}
\begin{proposition}\label{pro5}\cite{2005dsuhsVariational}
For any fixed $\varepsilon>0$,~$G_{\varepsilon}^{*}$ is the Fenchel-Legendre transformation of the soft-max function $G_{\varepsilon}$, then
\begin{equation}\label{a2.2.15}
\begin{array}{lll}
G_{\varepsilon}^{*}(\mathbf{w})&=\max\limits_{\mathbf{x}}\{<\mathbf{x},\mathbf{w}>-
G_{\varepsilon}(\mathbf{x})\}\\
&=\left\{
\begin{array}{lll}
\varepsilon\displaystyle\sum\limits_{l=1}^{N}w_{l}\log w_{l},& \mathbf{w}=(w_{1},w_{2},...,w_M)\in\Delta^+,\\
+\infty,& else,
\end{array}
\right.
\end{array}
\end{equation}
where $\Delta^+=\{\mathbf{w}=(w_{1},w_{2},.....,w_M)\mid0\leq w_{l}\leq 1,\sum\limits_{l=1}^{M}w_{l}=1\},$
and thus
\begin{equation}\label{dual_formula}
G_{\varepsilon}(\mathbf{x})=G_{\varepsilon}^{**}(\mathbf{x})=\max\limits_{\mathbf{w}\in\Delta^+}
\left\{<\mathbf{w},\mathbf{x}>-\varepsilon\displaystyle\sum\limits_{l=1}^{N}w_{l}\log w_{l}\right\},
\end{equation}
where $G_{\varepsilon}^{**}$ is the Fenchel-Legendre transformation of $G_{\varepsilon}^{*}$.
\end{proposition}
\par
Set $x_l=\varepsilon\log P_{\pmb{ \mathscr{V}}}(\pmb{f}_i-R_l\pmb{u}),~l=1,2,...,N$. For any fixed $i$ index, by Proposition \ref{pro5}, we can get
\begin{equation}\label{a2.2.16}
\begin{array}{lll}
   \log \displaystyle\sum_{l=1}^NP_{\pmb{ \mathscr{V}}}(\pmb{f}_i-R_l\pmb{u})&=& \frac{1}{\varepsilon}\cdot\varepsilon
  \log\displaystyle \sum_{l=1}^Ne^{\frac{\pmb{x}_l}{\varepsilon}}\\
   &=& \frac{1}{\varepsilon}\max\limits_{\mathbf{w}_i\in\Delta^+}
\left\{<\mathbf{w},\mathbf{x}>-\varepsilon\displaystyle\sum\limits_{l=1}^{N}w_{i,l}\log w_{i,l}\right\}\\
&=&\max\limits_{\mathbf{w}_i\in\Delta^+}\left\{\displaystyle\sum_{l=1}^Nw_{i,l}\log P_{\pmb{ \mathscr{V}}}(\pmb{f}_i-R_l\pmb{u})-\displaystyle\sum\limits_{l=1}^{N}w_{i,l}\log w_{i,l}\right\}.
\end{array}
\end{equation}
Substitute \eqref{a2.2.16} into \eqref{a2.2.12}, then we can get
\begin{equation}\label{a2.2.17}
\displaystyle\min_{\pmb{u},\pmb{\theta}}\min\limits_{\mathbf{w}_i\in\Delta^+}
\displaystyle\Bigg\{\mathcal{H}(\pmb{u},\pmb{\theta},\pmb{w})=-\sum_{i=1}^M\sum_{l=1}^N
w_{i,l}\log P_{\pmb{ \mathscr{V}}}(\pmb{f}_i-R_l\pmb{u})+\displaystyle\sum_{i=1}^M\sum\limits_{l=1}^{N}w_{i,l}\log w_{i,l}-\log P_{\pmb{\mathscr{U}}}(\pmb{u})\Bigg\}.
\end{equation}
In addition, we have
\begin{equation}\label{a2.2.18}
P_{\pmb{ \mathscr{V}}}(\pmb{f}_i-R_l\pmb{u})=\prod_{j=1}^NP_{ \mathscr{V}}(f_{i,j}-R_lu_j),~~~~P_{\pmb{\mathscr{U}}}(\pmb{u})=\prod_{c\in\mathcal{C}}p_{\mathscr{U}}(u_c),
\end{equation}
where $\mathcal{C}$ is the clicks in the graph representation of the prior.
Substitute \eqref{a1.6}, \eqref{a2.2.4}, \eqref{a2.2.18} into \eqref{a2.2.17} successively, we can get
\begin{equation}\label{a2.2.19}
\begin{array}{lll}
  & \displaystyle\min_{\pmb{u},\pmb{\theta}}\min\limits_{\mathbf{w}_i\in\Delta^+}
&\Bigg\{\mathcal{H}(\pmb{u},\pmb{\theta},\pmb{w})=-\displaystyle\sum_{i=1}^M
\sum_{l=1}^N
w_{i,l}\sum_{j=1}^N\log \sum_{k=1}^K\alpha_k p_{\mathscr{V}_k}(f_{i,j}-R_lu_j)
\\
 & &
 \displaystyle\sum_{i=1}^M\sum\limits_{l=1}^{N}w_{i,l}\log w_{i,l}+\gamma\sum_{c\in\mathcal{C}}\phi(u_c)\Bigg\}.
  \end{array}
\end{equation}
Because there is a $\log\sum$ item in \eqref{a2.2.19}, we apply Proposition \ref{pro5} again. Here let $N=K,$~$x_k=\varepsilon\log[\alpha_k p_{\mathscr{V}_k}(f_{i,j}-R_lu_j)]$ for any $k=1,2,...,K$. For fixed indexes $i,~j,~l$, we can get
\begin{equation}\label{a2.2.20}
\begin{array}{lll}
&&\log\displaystyle\sum_{k=1}^K\alpha_k p_{\mathscr{V}_k}(f_{i,j}-R_lu_j)
\\
&=&\displaystyle\max_{\pmb{q}_{i,j,l}\in Q^+}\Bigg\{
\sum_{k=1}^K\left[q_{i,j,l,k}(\log\alpha_k+\log p_{\mathscr{V}_k}(f_{i,j}-R_lu_j))\right]-\displaystyle\sum_{k=1}^Kq_{i,j,l,k}\log q_{i,j,l,k}\Bigg\},
\end{array}
\end{equation}
where $Q^+=\{\pmb{q}=(q_1,q_2,...,q_K)\mid0\leq q_k\leq 1,\displaystyle\sum_{k=1}^Kq_k=1\}$.
Then we substitute \eqref{a2.2.20} to\eqref{a2.2.19}. By simple calculation, we can get
\begin{equation}\label{a2.2.21.a1}
\begin{aligned}
\min_{\pmb{u},\pmb{\theta},\pmb{w}_{i}\in\Delta^+,\pmb{q}_{i,j,l}\in Q^+}
&\Bigg\{\mathcal{J}(\pmb{u},\pmb{\theta},\pmb{w},\pmb{q})=
-\sum_{i=1}^M\sum_{l=1}^Nw_{i,l}\sum_{j=1}^N\sum_{k=1}^K\left[q_{i,j,l,k}(\log\alpha_k+\log p_{\mathscr{V}_k}(f_{i,j}-R_lu_j))\right]\\
&+\sum_{i=1}^M\sum_{l=1}^Nw_{i,l}\sum_{j=1}^N\sum_{k=1}^Kq_{i,j,l,k}\log q_{i,j,l,k}+\displaystyle\sum_{i=1}^M\sum\limits_{l=1}^{N}w_{i,l}\log w_{i,l}
+\gamma\sum_{c\in\mathcal{C}}\phi(u_c)\Bigg\},
\end{aligned}
\end{equation}
where $\Delta^+=\{\bm{w}=(w_1,w_2,...,w_M)\mid 0\leq w_i\leq1,\displaystyle\sum_{i=1}^Mw_i=1\}$ and $Q^+=\{\pmb{q}=(q_1,q_2,...,q_K)\mid0\leq q_k\leq1,\displaystyle\sum_{k=1}^Kq_k=1\}$.
\subsubsection{Mixed Gaussian-Gaussian noise model}
Here we set $K=2$, $\mathscr{V}_1,~\mathscr{V}_2$ follow Gaussian distributions with mean $0$
and variances $\sigma_1^2,~\sigma_2^2$ respectively. The probability density function of $\mathscr{V}_k$ is
\begin{equation}\label{a1.3}
p_{\mathscr{V}_k}(v)=\frac{1}{\sqrt{2\pi \sigma_k^2}}e^{-\frac{v^2}{2\sigma_k^2}},
\end{equation}
for any $k=1,2$.
\par
Substitute \eqref{a1.3} into \eqref{a2.2.21.a1}. By simple calculation, we can get
\begin{equation}\label{a2.2.21}
\begin{aligned}
\min_{\pmb{u},\pmb{\theta},\pmb{w}_{i}\in\Delta^+,\pmb{q}_{i,j,l,k}\in Q^+}
&\Bigg\{\mathcal{J}(\pmb{u},\pmb{\theta},\pmb{w},\pmb{q})=
\sum_{i=1}^M\sum_{l=1}^Nw_{i,l}\Bigg[\sum_{j=1}^N\sum_{k=1}^2\frac{
(f_{i,j}-R_lu_j)^2}{2\sigma_k^2}q_{i,j,l,k}\\
&+\sum_{j=1}^N\sum_{k=1}^2\Big(\frac{1}{2}\log\sigma_k^2-\log\alpha_k\Big)
q_{i,j,l,k}+\sum_{j=1}^N\sum_{k=1}^2q_{i,j,l,k}\log q_{i,j,l,k}\Bigg]\\
&+\displaystyle\sum_{i=1}^M\sum\limits_{l=1}^{N}w_{i,l}\log w_{i,l}+\gamma\sum_{c\in\mathcal{C}}\phi(u_c)\Bigg\}.
\end{aligned}
\end{equation}
where $\Delta^+=\{\bm{w}=(w_1,w_2,...,w_M)\mid 0\leq w_i\leq1,\displaystyle\sum_{i=1}^Mw_i=1\}$ and $Q^+=\{\pmb{q}=(q_1,q_2)\mid0\leq q_k\leq1,\displaystyle\sum_{k=1}^2q_k=1\}$.
\section{Statistical interpretation of the proposed model}\label{Sec3}
The observed data is $A=(\pmb{f}_1,\pmb{f}_2,...,\pmb{f}_M)$. Recall that in Subsection \ref{subsec_modeling}, by maximum a posteriori(MAP) estimation, the goal becomes
\begin{equation}\label{sta_1}
\max_{\pmb{u},\pmb{\theta}}\left\{\mathcal{L}(\pmb{u},\pmb{\theta})\right\},
\end{equation}
where
\begin{equation}\label{sta_1_1}
 \mathcal{L}(\pmb{u},\pmb{\theta})=\log p(A|\pmb{u};\pmb{\theta})+\log p_{\pmb{\mathscr{U}}}(\pmb{u}).
\end{equation}
\par
Next, we will apply Expectation-Maximization(EM) algorithm to \eqref{sta_1}. Because $A$ is an incomplete data, we introduce a hidden variable $\pmb{y}=(\mathscr{L},\mathscr{K})$ to get the complete data $\pmb{Z}=(A,\pmb{y})$, where $\mathscr{L}$ and $\mathscr{K}$ represent the random variables of circularly shifted operator and the type of noise separately. $l$ and $c$ are some realizations of the random variables $\mathscr{L}$ and $\mathscr{K}$ separately. By Bayes' law, we get
\begin{equation}\label{sta_2}
p(\pmb{Z}|\pmb{u};\pmb{\theta})=p(A,\pmb{y}|\pmb{u};\pmb{\theta})=p\big(\pmb{y}|(A,\pmb{u});\pmb{\theta}\big)
p(A|\pmb{u};\pmb{\theta}).
\end{equation}
By taking the logarithm of \eqref{sta_2}, then
\begin{equation}\label{sta_3}
\log p(A|\pmb{u};\pmb{\theta})=\log p(\pmb{Z}|\pmb{u};\pmb{\theta})-\log p\big(\pmb{y}|(A,\pmb{u});\pmb{\theta}\big).
\end{equation}
Take the expectation of $\log p(A|\pmb{u};\pmb{\theta})$ with respect of $\pmb{Z}$ under the condition of $\pmb{\theta}^{\nu}$, where $\nu$ is the iteration step. We can get
\begin{equation}\label{sta_4}
E_{\pmb{y}}[\log p(\pmb{Z}|\pmb{u};\pmb{\theta})|(A,\pmb{u}^\nu);\pmb{\theta}^\nu]=\sum_{\pmb{y}}p(\pmb{y}|(A,\pmb{u}^\nu);
\pmb{\theta}^\nu)\log p(\pmb{Z}|\pmb{u};\pmb{\theta}),
\end{equation}
which is called E-step. Then we maximize the above expectation to get $\pmb{u}^{\nu+1},\pmb{\theta}^{\nu+1}$ as follows
\begin{equation}\label{sta_5}
(\pmb{u}^{\nu+1},\pmb{\theta}^{\nu+1})=\displaystyle\arg\max_{\pmb{u},\pmb{\theta}}\left\{E_{\pmb{y}}[\log p(\pmb{Z}|\pmb{u};\pmb{\theta})|(A,\pmb{u}^\nu);\pmb{\theta}^\nu]+\log p_{\pmb{\mathscr{U}}}(\pmb{u})\right\},
\end{equation}
which is called M-step.
\par
Calculating expectations is a key step in EM algorithm. Next, we will derive the concrete expression of \eqref{sta_4}.
\par
Firstly, we consider $\log p(\pmb{Z}|\pmb{u};\pmb{\theta})$. For any $\pmb{y}=(l,k)$, 
\begin{equation}\label{sta_7}
\begin{aligned}
\log p(\pmb{Z}|\pmb{u};\pmb{\theta})
=&\log\displaystyle\prod_{i=1}^Mp(\pmb{f}_i,\pmb{y}|\pmb{u};\pmb{\theta})
=\sum_{i=1}^M\log p(\pmb{f}_i,\pmb{y}|\pmb{u};\pmb{\theta})=\sum_{i=1}^M\log\prod_{j=1}^N p(f_{i,j},\pmb{y}|u_j;\pmb{\theta})\\
=&\sum_{i=1}^M\sum_{j=1}^N\log p(f_{i,j},\pmb{y}|u_j;\pmb{\theta})=\sum_{i=1}^M\sum_{j=1}^N
\left[\log p(\pmb{y}|u_j;\pmb{\theta})+\log p(f_{i,j}|(\pmb{y},u_j);\pmb{\theta})\right]\\
=&\sum_{i=1}^M
\sum_{j=1}^N\left[\log\left(\frac{1}{N}\alpha_k\right)+\log p_{\mathscr{V}_{k}}(f_{i,j}-R_{l}u_j;\pmb{\theta})\right],
\end{aligned}
\end{equation}
where $p\big((l,k)|u_j;\pmb{\theta}\big)=p(l,k)=p_{\mathscr{L}}(l))p_{\mathscr{C}}(k)=\frac{1}{N}\alpha_k$ owing to the independence of $\mathscr{L}$ and $\mathscr{K}$, and $p(f_{i,j}|(\pmb{y},u_j);\pmb{\theta})=p(f_{i,j}|(l,k,u_j);\pmb{\theta})=
p_{\mathscr{V}_{k}}(f_{i,j}-R_{l}u_j;\pmb{\theta}))$ for any $\pmb{y}=(l,k)$.
\par
Then we consider $p(\pmb{y}|(A,\pmb{u}^\nu);\pmb{\theta}^\nu)$. By Bayes' law, for any $\pmb{y}=(l,k)$, we can get
\begin{equation}\label{sta_10}
p(l,k|(A,\pmb{u}^\nu);\pmb{\theta}^\nu)=p(l|(A,\pmb{u}^\nu);\pmb{\theta}^\nu)p(k|(l,A,\pmb{u}^\nu);\pmb{\theta}^\nu).
\end{equation}
Specifically,
\begin{equation}\label{sta_11}
\begin{aligned}
p(l|(\pmb{f}_{i},\pmb{u}^\nu);\pmb{\theta}^\nu)=&\frac{p(l,\pmb{f}_{i}|\pmb{u}^\nu;
\pmb{\theta}^\nu)}{p(\pmb{f}_{i}|\pmb{u}^\nu;\pmb{\theta}^\nu)}
=\frac{p(\pmb{f}_{i}|(l,\pmb{u}^\nu);\pmb{\theta}^\nu)p_{\mathscr{L}}(l)}
{\displaystyle\sum_{l'=1}^Np(\pmb{f}_{i}|(l',\pmb{u}^\nu);\pmb{\theta}^\nu)p_{\mathscr{L}}(l')}\\
=&
\frac{p_{\pmb{\mathscr{V}}}(\pmb{f}_i-R_l\pmb{u}^\nu;\pmb{\theta}^\nu)}
{\displaystyle\sum_{l'=1}^Np_{\pmb{\mathscr{V}}}(\pmb{f}_i-R_{l'}\pmb{u}^\nu;\pmb{\theta}^\nu)}
=\frac{\displaystyle\prod_{j=1}^Np_{\mathscr{V}}(f_{i,j}-R_lu_j^\nu;\theta^\nu)}
{\displaystyle\sum_{l'=1}^N\displaystyle\prod_{j=1}^Np_{\mathscr{V}}(f_{i,j}-R_{l'}u_j^\nu;\theta^\nu)}
\doteq w_{i,l}^{\nu},
\end{aligned}
\end{equation}
where $p_{\mathscr{L}}(l)=\displaystyle\frac{1}{N}$ and
\begin{equation}\label{sta_11.1}
p_{\mathscr{V}}(f_{i,j}-R_{l}u_j^\nu;\theta^\nu)=\sum_{k=1}^K\alpha_k p_{\mathscr{V}_k}(f_{i,j}-R_{l}u_j^\nu;\theta^\nu)=\sum_{k=1}^K \frac{\alpha_k}{\sqrt{2\pi\sigma_k^{2(\nu)}}}\exp{\left(\frac{(f_{i,j}-R_{l}u_j^\nu)^2}
{2\sigma_k^{2(\nu)}}\right)}.
\end{equation}
\begin{equation}\label{sta_12}
\begin{aligned}
p(k|(l,f_{i,j},u_j^\nu);\pmb{\theta}^\nu)=&\frac{p\big(k,f_{i,j}|(l,u_j^\nu);\pmb{\theta}^\nu\big)}
{p\big(f_{i,j}|(l,u_j^\nu);\pmb{\theta}^\nu\big)}
=\frac{p\big(f_{i,j}|(k,l,u_j^\nu);\pmb{\theta}^\nu\big)p\big(k|(l,u_j^\nu);\pmb{\theta}^\nu\big)}
{\displaystyle\sum_{k'=1}^K
p\big(f_{i,j}|(k',l,u_j^\nu);\pmb{\theta}^\nu\big)p\big(k'|(l,u_j^\nu);\pmb{\theta}^\nu\big)}\\
=&\frac{\alpha_kp_{\mathscr{V}_k}(f_{i,j}-R_lu_j^\nu)}{\displaystyle\sum_{k'=1}^K
\alpha_{k'}p_{\mathscr{V}_{k'}}(f_{i,j}-R_lu_j^\nu)}\doteq q_{i,j,l,k}^\nu,
\end{aligned}
\end{equation}
where $p\big(f_{i,j}|(k,l,u_j^\nu);\pmb{\theta}^\nu\big)=p_{\mathscr{V}_k}(f_{i,j}-R_lu_j^\nu)$ and $p\big(k|(l,u_j^\nu);\pmb{\theta}^\nu\big)=\alpha_k$. Substitute \eqref{sta_11} and\eqref{sta_12} into \eqref{sta_10}, then we can get
\begin{equation}\label{sta_13}
p(l,k|f_{i,j};\pmb{\theta}^\nu)=w_{i,l}^\nu q_{i,j,l,k}^\nu.
\end{equation}
Substitute \eqref{sta_7} and \eqref{sta_13} into \eqref{sta_4}, the expectation \eqref{sta_4} becomes
\begin{equation}\label{sta_14}
E_{\pmb{y}}[\log p(\pmb{Z}|\pmb{u};\pmb{\theta})|(A,\pmb{u}^\nu);\pmb{\theta}^\nu]
=\sum_{i=1}^M\sum_{j=1}^N\sum_{l=1}^N\sum_{k=1}^K
w_{i,l}^\nu q_{i,j,l,k}^\nu\log\left(\frac{1}{N}\alpha_{k}
p_{\mathscr{V}_{k}}\big(f_{i,j}-R_{l}u_j;\pmb{\theta})\right).
\end{equation}
Substitute $K=2$ and
\begin{equation}\label{sta_15}
p_{\mathscr{V}_{k}}(f_{i,j}-R_{l}u_j;\pmb{\theta}))=\frac{1}{\sqrt{2\pi\sigma_k^2}}\exp
{\left(-\frac{\|f_{i,j}-R_lu_j\|_2^2}{2\sigma_k^2}\right)}
\end{equation}
into \eqref{sta_14},then we can get the iterative formula as follows
\begin{equation}\label{sta_16}
\begin{aligned}
&E_{\pmb{y}}[\log p(\pmb{Z}|\pmb{u};\pmb{\theta})|(A,\pmb{u}^\nu);\pmb{\theta}^\nu]+\gamma\sum_{c\in\mathcal{C}}\phi(u_c)\\
=&\sum_{i=1}^M\sum_{j=1}^N\sum_{l=1}^N\sum_{k=1}^2
w_{i,l}^\nu q_{i,j,l,k}^\nu\left(\log\alpha_{k}-\frac{1}{2}
\log\sigma_k^2-\frac{\|f_{i,j}-R_lu_j\|_2^2}{2\sigma_k^2}\right)\\
&+\gamma\sum_{c\in\mathcal{C}}\phi(u_c)\doteq\mathcal{T}(\pmb{u},\pmb{\theta}|\pmb{u}^\nu;\pmb{\theta}^\nu).
\end{aligned}
\end{equation}
Moreover, the M-step is
\begin{equation}\label{sta_17}
(\pmb{u}^{\nu+1},\pmb{\theta}^{\nu+1})=\arg\max_{\pmb{u},\pmb{\theta}}\mathcal{T}(\pmb{u},\pmb{\theta}|\pmb{u}^\nu
;\pmb{\theta}^\nu).
\end{equation}
\par
\section{The existence of a minimizer}\label{Sec4}
In this section, we prove the existence of a minimizer for the proposed model with total variation(TV) regularizer.

For the convenience of discussion, we introduction the following denotations. Denote $\Omega\subset\mathbb{R}$ as a bound open sets. $\mathbb{T},\mathbb{X}\in\mathbb{R}^+$ are bounded close sets. $BV(\Omega)$ is denoted as the space of bounded total variation function, i.e., \\
$BV(\Omega)=\{u\in L^1(\Omega)|J(u)<+\infty\},$ where
$J(u)=\sup\{\int_{\Omega}u(x)div(\varphi(x))dx|\varphi\in C_0^\infty(\Omega,\mathbb{R}^2), $ $\|\varphi\|_{L^\infty(\Omega,\mathbb{R}^2)}\leq1\}.$
\\
$S(\Omega)=\{u\in BV(\Omega)|u\geq0\}$.\\
$\mathbb{K}=\{\pmb{\theta}(t)=(\sigma^2(t),\alpha(t))|0
\leq\sigma_{min}^2\leq\sigma^2(t)\leq\sigma_{max}^2,0\leq\alpha_{min}\leq\alpha(t)<1, \forall t\in\mathbb{T}\}$.\\
$\mathbb{W}=\{w(x,y)\in L^\infty(\mathbb{X}\times\Omega)|0\leq w(x,y)dy\leq 1,\forall (x,y)\in\mathbb{X}\times\Omega;\int_{y}w(x,y)=1,\forall x\in\mathbb{X}\}.$
\\
$\mathbb{Q}=\{q(x,y,s,t)|0\leq q(x,y,s,t)\leq 1,\forall (x,y,s,t)\in\mathbb{X}\times\Omega\times\Omega\times\mathbb{T};\int_{t}q(x,y,s,t)=1,\forall (x,y,s)\in\mathbb{X}\times\Omega\times\Omega\}.$\\
$\Gamma=\{(u,\pmb{\theta},w,q)|u\in S(\Omega),\pmb{\theta}\in\mathbb{K},w\in\mathbb{W},q\in\mathbb{Q}\}$.
Denote $\hat{u}$ as an extension of $u$, i.e., $\hat{u}(z)= u(z)$, if~$ z\in\Omega$;  $\hat{u}(z)=u(z+N)$, if~$z\in\Omega-N=\{x-N|x\in\Omega\}$. The continuous form of the proposed model can be written as
\\
 \begin{equation}\label{b2.2.21}
\begin{aligned}
&\min_{(u,\pmb{\theta},w,q)}
\Bigg\{\mathcal{J}(u,\pmb{\theta},w,q)=
\int_{\mathbb{X}}\int_{\Omega}w(x,y)\Big[\int_{\Omega}\int_{\mathbb{T}}\frac{
(f(x,s)-\hat{u}(s-y))^2}{2\sigma(t)^2}
q(x,y,s,t)dtds\\
&+\int_{\Omega}\int_{\mathbb{T}}\Big(\frac{1}{2}\log\sigma^2(t)-\log\alpha(t)\Big)
q(x,y,s,t)dtds+\int_{\Omega}\int_{\mathbb{T}}q(x,y,s,t)\log q(x,y,s,t)dtds\Big]\\
&dydx+\int_{\mathbb{X}}\int_{\Omega}w(x,y)\log w(x,y)dydx+\lambda \phi(u)\Bigg\}.
\end{aligned}
\end{equation}
\\
\begin{theorem}\label{theorem1}
Assume~$f\in L^\infty(\mathbb{X}\times\Omega),~\displaystyle0<\inf_{\mathbb{X}\times\Omega}f,\sup_{\mathbb{X}\times\Omega}f<+\infty$. Set $\phi(u)=J(u)$; then there exists at least a solution in $\Gamma$ for the problem \eqref{b2.2.21}.
\end{theorem}
\par
The proof is in Appendix A.
\section{Related algorithm}\label{Sec5}
In this section, we design an algorithm for problem \eqref{a2.2.21}. Recall that problem \eqref{a2.2.21} is written as
\begin{equation*}
\begin{aligned}
\min_{\pmb{u},\pmb{\theta},\pmb{w}_{i}\in\Delta^+,\pmb{q}_{i,j,l}\in Q^+}
&\Bigg\{\mathcal{J}(\pmb{u},\pmb{\theta},\pmb{w},\pmb{q})=
\sum_{i=1}^M\sum_{l=1}^Nw_{i,l}\Bigg[\sum_{j=1}^N\sum_{k=1}^2\frac{
(f_{i,j}-R_lu_j)^2}{2\sigma_k^2}q_{i,j,l,k}\\
&+\sum_{j=1}^N\sum_{k=1}^2\Big(\frac{1}{2}\log\sigma_k^2-\log\alpha_k\Big)
q_{i,j,l,k}+\sum_{j=1}^N\sum_{k=1}^2q_{i,j,l,k}\log q_{i,j,l,k}\Bigg]\\
&+\displaystyle\sum_{i=1}^M\sum\limits_{l=1}^{N}w_{i,l}\log w_{i,l}+\gamma\sum_{c\in\mathcal{C}}\phi(u_c)\Bigg\}.
\end{aligned}
\end{equation*}
 Denote $\nu$ as outer iteration steps. By the alternating direction iterative method, we can decompose the problem \eqref{a2.2.21} into two subproblems as follows,
\begin{equation}\label{a5.1.1}
 \left \{\begin{array}{l}
 (\pmb{u}^{\nu+1},\pmb{\theta}^{\nu+1})=\displaystyle\arg\min_{\pmb{u},\pmb{\theta}}\mathcal{J}(\pmb{u},
 \pmb{\theta},\pmb{w}^{\nu},\pmb{q}^\nu),\\
  (\pmb{w}_i^{\nu+1},\pmb{q}_{i,j,l}^{\nu+1})=\displaystyle\arg\min_{\pmb{w}_i,\pmb{q}_{i,j,l}}
 \mathcal{J}(\pmb{u}^{\nu+1},\pmb{\theta}^{\nu+1},\pmb{w},\pmb{q}).
  \\
 \end{array}\right.
\end{equation}
We can get a close solution of $\pmb{w}_i^{\nu+1},$ by taking the derivative of $\mathcal{J}(\pmb{u}^{\nu+1},\pmb{\theta}^{\nu+1},\pmb{w},\pmb{q})$ with respect of $\pmb{w}_i$ as follows
\begin{equation}\label{a5.1.2}
w_{i,l}^{\nu+1}=\frac{p_{\pmb{\mathscr{V}}}(\pmb{f}_i-R_l\pmb{u}^{\nu+1})}
{\displaystyle\sum_{l'=1}^N
p_{\pmb{\mathscr{V}}}(\pmb{f}_i-R_{l'}\pmb{u}^{\nu+1})},~~~~~l=1,2,...,N,
\end{equation}
where
\begin{equation}\label{a5.1.3}
\begin{array}{ll}
 p_{\pmb{\mathscr{V}}}(\pmb{f}_i-R_l\pmb{u}^{\nu+1})&=\displaystyle\prod_{j=1}^Np_{\mathscr{V}}(f_i-R_lu^{\nu+1})
 \\
 &=\displaystyle\prod_{j=1}^N\sum_{k=1}^2\alpha_k^{\nu+1} p_{\mathscr{V}_k}(f_{i,j}-R_lu_j^{\nu+1})\\
 &=\displaystyle\prod_{j=1}^N\left[\sum_{k=1}^2 \frac{\alpha_k^{\nu+1}
 }{\sqrt{2\pi\sigma_k^{2(\nu+1)}}}\exp{\left(-\frac{(f_{i,j}-R_lu_j^{\nu+1})^2}{2
 \sigma_k^{2(\nu+1)}}\right)}\right]\\
 \end{array}
\end{equation}
We can get a close solution of $q_{i,j,l,k}^{\nu}$ by taking the derivative of $\mathcal{J}(\pmb{u}^{\nu+1},\pmb{\theta}^{\nu+1},\pmb{w},\pmb{q})$  with the respect of $q_{i,j,l,k}$ as follows,
\begin{equation}\label{a5.1.8}
\begin{array}{ll}
  q_{i,j,l}^{\nu+1}\doteq q_{i,j,l,1}^{\nu+1}&=\displaystyle\frac{\alpha_k^{\nu+1} p_{\mathscr{V}_k}(f_{i,j}-R_lu_j^{\nu+1})}{\displaystyle\sum_{k'=1}^K\alpha_{k'}^{\nu+1} p_{\mathscr{V}_{k'}}(f_{i,j}-R_lu_j^{\nu+1})},\\
&=\displaystyle\frac{ \frac{\alpha_k^{\nu+1}}{\sigma_k^{\nu+1}}
\exp{\left(-\frac{\left(f_{i,j}-R_lu_j^{\nu+1}\right)^2}
{2\sigma_1^{2(\nu+1)}}\right)}}
{\displaystyle\sum_{k'=1}^2\frac{\alpha_{k'}^{\nu+1}}{\sigma_{k'}^{\nu+1}}
\exp{\left(-\frac{\left(f_{i,j}-R_lu_j^{\nu+1}\right)^2}
{2\sigma_{k'}^{2(\nu+1)}}
\right)}},
\end{array}
\end{equation}
where $i=1,2,...,M;~j,l=1,2,...,N$.
\par
 It's easy to see that the first subproblem of \eqref{a5.1.1} equals to the problem \eqref{sta_17}. In other wards, the iterative scheme \eqref{a5.1.1} derived by soft-max method is same as the iterative scheme \eqref{sta_17} derived by EM algorithm.
\par
Now we come back to the first subproblem of \eqref{a5.1.1}. The discussion below is in the case of fixed outer iteration steps $\nu$.  By the alternating direction iterative method, we can get
\begin{equation}\label{a5.1.9}
\left\{\begin{array}{l}
\pmb{u}^{\nu}=\displaystyle\arg\min_{\pmb{u}}
  \mathcal{J}(\pmb{u},\pmb{\theta}^{\nu},\pmb{w}^{\nu},\pmb{q}^{\nu}),\\
\pmb{\theta}^{\nu+1}=\displaystyle\arg\min_{\pmb{u}}
  \mathcal{J}(\pmb{u}^{\nu+1},\pmb{\theta},\pmb{w}^{\nu},\pmb{q}^{\nu}).
\end{array}\right.
\end{equation}
We can get the following iterative formulas from the second subproblem of \eqref{a5.1.9}:
\begin{equation}\label{a5.1.11}
\alpha_k^{\nu+1}=\frac{\displaystyle\sum_{i=1}^M\sum_{l=1}^Nw_{i,l}^{\nu}\left(\sum_{j=1}^N
{q_{i,j,l,k}^{\nu}}\right)}
{MN},
\end{equation}
\begin{equation}\label{a5.1.12}
\sigma_k^{2(\nu+1)}=\frac{\displaystyle\sum_{i=1}^M\sum_{l=1}^Nw_{i,l}^{\nu}\left(\sum_{j=1}^N
\left(f_{i,j}-R_lu_j^{\nu+1}\right)^2q_{i,j,l,k}^{\nu}\right)}
{\displaystyle\sum_{i=1}^M\sum_{l=1}^Nw_{i,l}^{\nu}\sum_{j=1}^N\sum_{k'=1}^2q_{i,j,l,k'}^{\nu}},
\end{equation}
\par
For the first subproblem of \eqref{a5.1.9}, we can rewrite it as
\begin{equation}\label{a5.1.16}
\begin{array}{ll}
\pmb{u}^{\nu+1}=\displaystyle\arg\min_{\pmb{u}}
&\Bigg\{\displaystyle\sum_{i=1}^M\sum_{l=1}^Nw_{i,l}^{\nu}\sum_{j=1}^N\sum_{k=1}^2\frac{
\left(f_{i,j}-R_lu_j\right)^2}{2(\sigma_k^2)^{\nu}}q_{i,j,l,k}^{\nu}
+\gamma \phi(\mathbf{u})\Bigg\},
\end{array}
\end{equation}
where $\phi(\mathbf{u})$ is a regularizer and $\gamma>0$ is a parameter. We apply the augmented Lagrange method to problem \eqref{a5.1.16} to get
\begin{equation}\label{a5.1.17a}
\begin{array}{ll}
\displaystyle\min_{\pmb{u},\pmb{d}}\max_{\pmb{p}}&\Bigg\{
\displaystyle\sum_{i=1}^M\sum_{l=1}^Nw_{i,l}^{\nu}\sum_{j=1}^N\sum_{k=1}^2\frac{
\left(f_{i,j}-R_ld_j\right)^2}{2(\sigma_k^2)^{\nu}}q_{i,j,l,k}^{\nu}\\
&+<\pmb{p},\pmb{u-d}>+\displaystyle\frac{r}{2}\|
\pmb{u}-\pmb{d}\|_2^2
+\gamma \phi(\mathbf{u})\Bigg\}.
\end{array}
\end{equation}
Fix the index $\nu$. Denote $\iota$ the inner iteration steps of subproblem \eqref{a5.1.17a}. The problem \eqref{a5.1.17a} can be decomposed into
\begin{equation}\label{a5.1.18}
\left\{\begin{array}{rl}
(\mathbf{u}^{\iota+1},\mathbf{d}^{\iota+1})=&\displaystyle\arg\min_{\mathbf{u},\mathbf{d}}\Bigg\{\displaystyle
  \sum_{i=1}^M\sum_{l=1}^Nw_{i,l}^{\nu}\sum_{j=1}^N\sum_{k=1}^2\frac{
\left(f_{i,j}-R_ld_j\right)^2}{2\sigma_k^{2(\nu)}}q_{i,j,l,k}^{\nu}\\
&+\displaystyle\frac{r}{2}\left\|
\mathbf{u}-\mathbf{d}+\frac{\mathbf{p}^{\iota}}{r}\right\|_2^2+\lambda \phi(\mathbf{u})\Bigg\},\\
  \mathbf{p}^{\iota+1}=&\mathbf{p}^{\iota}+\tau(\mathbf{u}^{\iota}-\mathbf{d}^{\iota}).
\end{array}\right.
\end{equation}
Furthermore, the problem \eqref{a5.1.18} can be decomposed once again and we can get
\begin{equation}\label{a5.1.19}
 \left\{\begin{array}{rl}
\mathbf{u}^{\iota+1}=&\displaystyle\arg\min_{\mathbf{u}}\Bigg\{\displaystyle\frac{r}{2}\left\|
\mathbf{u}-\mathbf{d}^\iota+\frac{\mathbf{p}^{\iota}}{r}\right\|_2^2+\lambda \phi(\mathbf{u})\Bigg\},\\
\mathbf{d}^{\iota+1}=&\displaystyle\arg\min_{\mathbf{d}}\Bigg\{\displaystyle
  \sum_{i=1}^M\sum_{l=1}^Nw_{i,l}^{\nu}\sum_{j=1}^N\sum_{k=1}^2\frac{
\left(f_{i,j}-R_ld_j\right)^2}{2\sigma_k^{2(\nu)}}q_{i,j,l,k}^{\nu}\\
&+\displaystyle\frac{r}{2}\left\|
\mathbf{d}-\mathbf{u}^{\iota+1}-\frac{\mathbf{p}^{\iota}}{r}\right\|_2^2\Bigg\},\\
  \mathbf{p}^{\iota+1}=&\mathbf{p}^{\iota}+\tau(\mathbf{u}^{\iota}-\mathbf{d}^{\iota}).
 \end{array}\right.
\end{equation}
The first subproblem of \eqref{a5.1.19} is a Gaussian denoiser. For the second subproblem of \eqref{a5.1.19}, we take the derivative with respect of $\pmb{d}$.
 We can get
\begin{equation}\label{a5.1.20}
d_j^{\iota+1}=\displaystyle\frac{\displaystyle\sum_{i=1}^M\sum_{l=1}^Nw_{i,l}^\nu t_1
+\sigma_1^{2(\nu)}\sigma_2^{2(\nu)}\left(ru_j^{\iota+1}+p_j^\iota\right)}
{\displaystyle\sum_{i=1}^M\sum_{l=1}^Nw_{i,l}^\nu t_2+r\sigma_1^{2(\nu)}
\sigma_2^{2(\nu)}},
\end{equation}
where
\begin{equation}\label{a5.1.21}
t_1=
\left[q_{i,j,l}^{\nu}\sigma_2^{2(\nu)}+
(1-q_{i,j,l}^{\nu})\sigma_1^{2(\nu)}\right]R_l^{-1}f_{i,j},
\end{equation}
\begin{equation}\label{a5.1.22}
t_2=
q_{i,j,l}^{\nu}\sigma_2^{2(\nu)}+(1-q_{i,j,l}^{\nu}).
\sigma_1^{2(\nu)}.
\end{equation}
According to the above discussion, We propose the following algorithm.
 \begin{algorithm}\label{alg1}
\textbf{Algorithm~1}~(MGG SoftMax)\\
\textbf{------------------------------------------------------------------------------------------------------------------------
}\\
1: Initialization. Let $\nu=0$. Set $\pmb{u}^0=\pmb{f}_1,\pmb{\theta}^0$. Then calculate $\pmb{w}^0$ by \eqref{a5.1.2} and calculate $\pmb{q}^0$ by \eqref{a5.1.8}.\\
 2:Smothness. Set the inner iteration $\iota=0, \pmb{u}^{\nu+1,0}=\pmb{u}^{\nu}$. Update $\pmb{u}^{\nu+1,\iota+1}$ by \eqref{a5.1.19} until convergence.  \\
3:Parameter estimation. Update $\alpha^{\nu+1}$ and  $\sigma_k^{2(\nu+1)}$
by\eqref{a5.1.11} and \eqref{a5.1.12} separately.\\
4:Noise classification. Update $\pmb{q}^{\nu+1}$ by calculating \eqref{a5.1.8}.\\
 5:Circularly-shifted classification. Update $\pmb{w}^{\nu+1}$ by calculating \eqref{a5.1.2}.\\
 6:convergence condition. If $\frac{\|\pmb{u}^{\nu+1}-\pmb{u}^{\nu}\|}{\|\pmb{u}^{\nu}\|}<\varepsilon$,  stop iterating the algorithm; Else, go to the step $2$. \\
 \end{algorithm}
\par
\section{Convergence analysis}\label{Sec6}
In this section, we will show some convergence analysis for the proposed algorithm.
\par
\begin{theorem}\label{theorem2}\cite{2013Adsiusuyy}
The functional $\mathcal{L}$,~$\mathcal{H}$ and $\mathcal{J}$ have the same minimizer $(\pmb{u}^*,\pmb{\theta}^*)$.~
\end{theorem}
\begin{theorem}\label{theorem3}\cite{2013Adsiusuyy}
(Energy Descent)~If the sequence $(\pmb u^\nu,\pmb \theta^\nu)$ satisfies $\mathcal{J}(\pmb u^{\nu+1},\pmb \theta^{\nu+1})\leq\mathcal{J}(\pmb u^{\nu},\pmb \theta^{\nu})$; then we can get
\begin{equation}\label{AB1}
\mathcal{L}(\pmb u^{\nu+1},\pmb \theta^{\nu+1})\leq\mathcal{L}(\pmb u^{\nu},\pmb \theta^{\nu}).
\end{equation}
\end{theorem}
\begin{theorem}\label{theorem5}\cite{2019Blockcsuih}
Fixed the iteration steps $\nu$.~Assume~ $\|\pmb{f}\|_{\infty}<+\infty~$. Set $\pmb u^*$~be the minimizer of the problem \eqref{a5.1.16}. $\forall~0< \tau< 2r$, then the sequence
~$\pmb u^{\iota}$~generated by the iteration scheme \eqref{a5.1.18} converges to $\pmb u^*$,~i.e.~$\displaystyle\lim_{\iota\rightarrow+\infty}\pmb u^{\iota}=\pmb u^*$.
\end{theorem}
\begin{theorem}\label{theorem6}\cite{2019Blockcsuih}
Fixed the iteration steps $\nu$.~Assume~ $\|\pmb{f}\|_{\infty}<+\infty~$. Set $\pmb u^*$~be the minimizer of the problem \eqref{a5.1.16}. let $\tau=r$, then the sequence
~$\pmb u^{\iota}$~generated by the iteration scheme \eqref{a5.1.19} converges to
$\pmb u^*$,~i.e.~$\displaystyle\lim_{\iota\rightarrow+\infty}\pmb u^{\iota}=\pmb u^*$.
\end{theorem}
\section{Numerical Experiments and Results}\label{Sec7}
This section is focused on numerical experiments. We consider the case of mixed Gaussian-Gaussian(MGG) noise. The true signal $\pmb{u}$ is a 1-dimensional signal of length $N$. We can generate $M$ noisy and circularly-shifted observations by the formula \eqref{a1.1}.
The goal is to recover the true signal $\pmb{u}$ from $M$ noisy circularly-shifted observations. All methods are evaluated by relative recovery error defined as
\begin{equation}\label{numer2}
  \text{Relative Error}(\tilde{\pmb{u}},\pmb{u})=\frac{\|R_l\tilde{\pmb{u}}-\pmb{u}\|_2}{\|\pmb{u}\|_2}
\end{equation}
where $\tilde{\pmb{u}}$ is an estimation of the true signal $\pmb{u}$.
\par
 To show the validity of the proposed method, we make a number of comparisons between the proposed model and some existing algorithms like Expectation Maximization(EM)\cite{katsevich2020likelihood}, spectral method largest spectral gap\cite{8352518}, frequency marching(FM)\cite{2017Bispectrum}, optimization on phase manifold\cite{2017Bispectrum} and optimization on phase synchronization\cite{2017Bispectrum}. The proposed model is abbreviated by MGG SoftMax. The experiments are run on a computer with $8$ Inter Core i7-8550U CPUs. These CPUs are used to compute thousands of FFTs for the proposed method and EM algorithm, while they are also used to compute the invariants in parallels.
\par
\begin{table}[!htb]
\renewcommand\arraystretch{1.25}
  \centering
  \caption{Comparison of relative recovery error values for a fixed real random signal $\pmb{u}$ of length $N=41$
  under different noise levels for EM\cite{katsevich2020likelihood},
  spectral M. largest spectral gap\cite{8352518},
  optim. phase manifold\cite{2017Bispectrum}, FM\cite{2017Bispectrum},
  iter. phase synch\cite{2017Bispectrum}
  and MGG SoftMax model with setting different initial values.
  (The largest relative errors are shown in bold fonts.)}
\label{TABLE1}
  \begin{tabular}{
c@{\hspace{3mm}}c@{\hspace{3mm}}c@{\hspace{3mm}}c@{\hspace{3mm}}
c@{\hspace{3mm}}c@{\hspace{3mm}}c@{\hspace{3mm}}c@{\hspace{3mm}}c@{\hspace{3mm}}
c@{\hspace{3mm}}c@{\hspace{3mm}}c@{\hspace{3mm}}c@{\hspace{3mm}}
c@{\hspace{3mm}}c@{\hspace{3mm}}c@{\hspace{3mm}}c@{\hspace{3mm}}
c c c c c c c c c c c c c c c c c}
    \hline
\hline
\\
\hline
$\alpha$&$\sigma_1$&$\sigma_2$&&
\multicolumn{3}{c}{Existing Methods} &&&
\multicolumn{1}{c}{Proposed Method}\\
\cline{5-11}
\Xcline{4-8}{1pt}\Xcline{10-11}{1pt}
\multirow{2}*{$\downarrow$}&\multirow{2}*{$\downarrow$}&\multirow{2}*{$\downarrow$}
&\multirow{2}*{EM\cite{katsevich2020likelihood}}
&\multicolumn{4}{c}{Invariant features}&&\multirow{2}*{MGG SoftMax}
\\
\cdashline{5-8}
&&&&Spec. M.\cite{8352518}&Opti. P. M.\cite{2017Bispectrum}&FM\cite{2017Bispectrum}
&Iter. P. S.\cite{2017Bispectrum}&&&
\\
\hline\hline
0&10&0.01
&0.0002082
&\pmb{0.0002078}
&0.0002079
&0.001146
&0.0002079
&&0.001267
\\
0.2&10&0.01
&0.7445
&0.6365
&0.6702
&0.6208
&0.6702
&&\pmb{0.08204}
\\
0.4&10&0.01
&0.7450
&0.7144
&0.6419
&0.6958
&0.6419
&&\pmb{0.1386}
\\
0.6&10&0.01
&0.7448
&0.7185
&0.7222
&0.7221
&0.7222
&&\pmb{0.2881}
\\
0.8&10&0.01
&0.7453
&0.7422
&0.7424
&0.7422
&0.7424
&&\pmb{0.6093}
\\
1&10&0.01
&0.7456
&0.7318
&0.7318
&0.7317
&0.7318
&&\pmb{0.6565}
\\
\hline
0&10&0.1
&\pmb{0.002082}
&0.003909
&0.002095
&0.01099
&0.002095
&&0.002758
\\
0.2&10&0.1
&0.7445
&0.6438
&0.6631
&0.6208
&0.6631
&&\pmb{0.08084}
\\
0.4&10&0.1
&0.745
&0.7136
&0.6552
&0.7056
&0.6552
&&\pmb{0.1668}
\\
0.6&10&0.1
&0.7448
&0.7186
&0.7224
&0.7208
&0.7224
&&\pmb{0.3355}
\\
0.8&10&0.1
&0.7453
&0.7422
&0.7423
 &0.7421
&0.7423
&&\pmb{0.5663}
\\
1&10&0.1
&0.7456
&0.7318
&0.7318
&0.7317
&0.7318
&&\pmb{0.6565}
\\
\hline
0&10&0.5
&\pmb{0.01046}
&0.1356
&0.01513
&0.03417
&0.03417
&&0.07679
\\
0.2&10&0.5
&0.7446
&0.6496
&0.6738
&0.6279
&0.6738
&&\pmb{0.09460}
\\
0.4&10&0.5
&0.7450
&0.6943
&0.6780
&0.7260
&0.6780
&&\pmb{0.5590}
\\
0.6&10&0.5
&0.7448
&0.7219
&0.7299
&0.7203
&0.7299
&&\pmb{0.5755}
\\
0.8&10&0.5
&0.7453
&0.7455
&0.7455
&0.7451
&0.7446
&&\pmb{0.6380}
\\
1&10&0.5
&0.7456
&0.7318
&0.7318
&0.7317
&0.7318
&&\pmb{0.6565}
\\
\hline
0&5&0.01
&0.0002082
&\pmb{0.0002078}
&0.0002079
&0.001146
&0.0002079
&&0.001267
\\
0.2&5&0.01
&0.7444
&0.6200
&0.5910
&0.6114
&0.5910
&&\pmb{0.05596}
\\
0.4&5&0.01
&0.7445
&0.5791
&0.6411
&0.6309
&0.6411
&&\pmb{0.07222}
\\
0.6&5&0.01
&0.7445
&0.6759
&0.7167
&0.6766
&0.7167
&&\pmb{0.2732}
\\
0.8&5&0.01
&0.7446
&0.6576
&0.6488
&0.6696
&0.6488
&&\pmb{0.5038}
\\
1&5&0.01
&0.7447
&0.6415
&\pmb{0.6366}
&0.7651
&0.6366
&&0.6524
\\
\hline
0&5&0.1
&\pmb{0.002082}
&0.003909
&0.002095
&0.01099
&0.002095
&&0.002758
\\
0.2&5&0.1
&0.7444
&0.6267
&0.6176
&0.6094
&0.6176
&&\pmb{0.05950}
\\
0.4&5&0.1
&0.7445
&0.6581
&0.6468
&0.6216
&0.6468
&&\pmb{0.08918}
\\
0.6&5&0.1
&0.7445
&0.6795
&0.7220
&0.5924
&0.7220
&&\pmb{0.2975}
\\
0.8&5&0.1
&0.7446
&0.6564
&0.6475
&0.6653
&0.6475
&&\pmb{0.5262}
\\
1&5&0.1
&0.7447
&0.6415
&\pmb{0.6366}
&0.7651
&0.6366
&&0.6524
\\
\hline
0&5&0.5
&\pmb{0.01046}
&0.1356
&0.01513
&0.03417
&0.01513
&&0.07679
\\
0.2&5&0.5
&0.7444
&0.6746
&0.5994
&0.6350
&0.5994
&&\pmb{0.07043}
\\
0.4&5&0.5
&0.7445
&0.6525
&0.5980
&0.6228
&0.5980
&&\pmb{0.5857}
\\
0.6&5&0.5
&0.7445
&0.6514
&0.6878
&0.5951
&0.6878
&&\pmb{0.5138}
\\
0.8&5&0.5
&0.7446
&0.6711
&0.6637
&0.6802
&0.6637
&&\pmb{0.6577}
\\
1&5&0.5
&0.7447
&0.6415
&\pmb{0.6366}
&0.7651
&0.6366
&&0.6524
\\
\hline
\end{tabular}
\end{table}
\begin{figure}[!htb]
\begin{minipage}[t]{0.49\linewidth}
\includegraphics[width=3.15in]{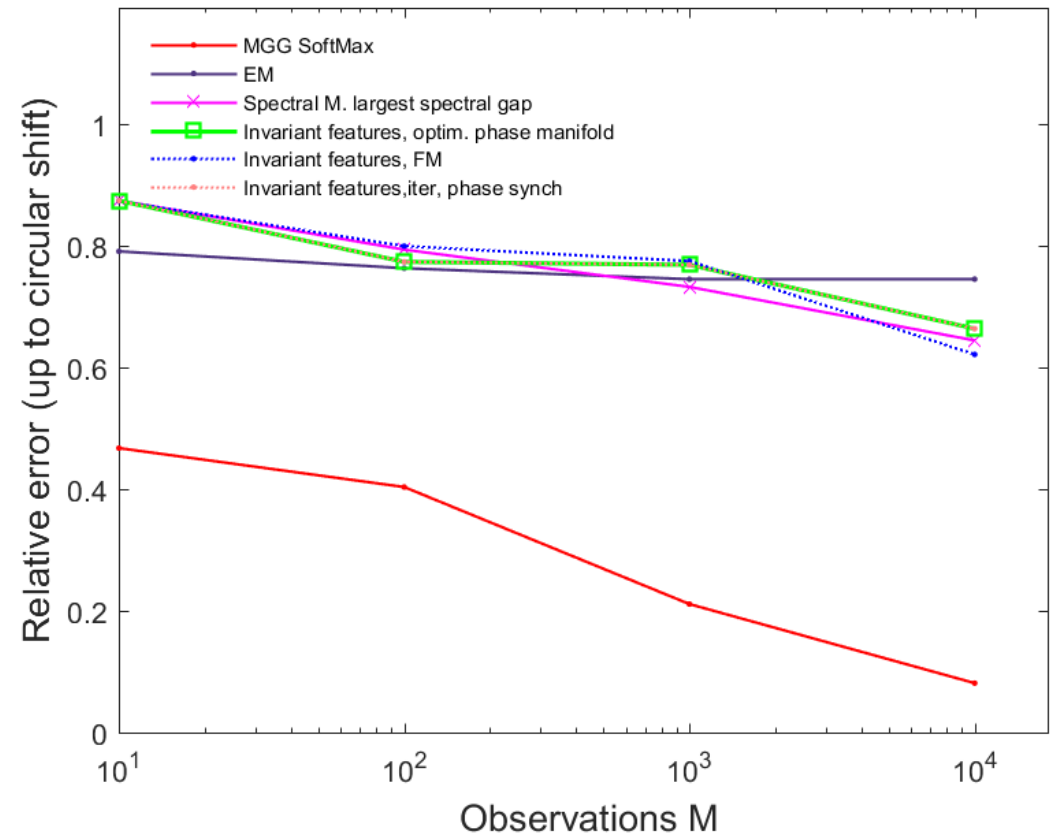}
\caption{Relative error as a function of the number of observations $M$. The observations are obtained through corrupting a real signal of length $N=41$ by a mixed Gaussian-Gaussian noise with fixed parameters $\alpha=0.2,\sigma_1=10,\sigma_2=0.1$.  Note that curves corresponding to the optim. phase manifold and the iter. phase synch almost overlap.}
\label{Fig1}
\end{minipage}%
\vspace{0.05 in}
\hspace{0.05 in}
\begin{minipage}[t]{0.49\linewidth}
\includegraphics[width=3.15in]{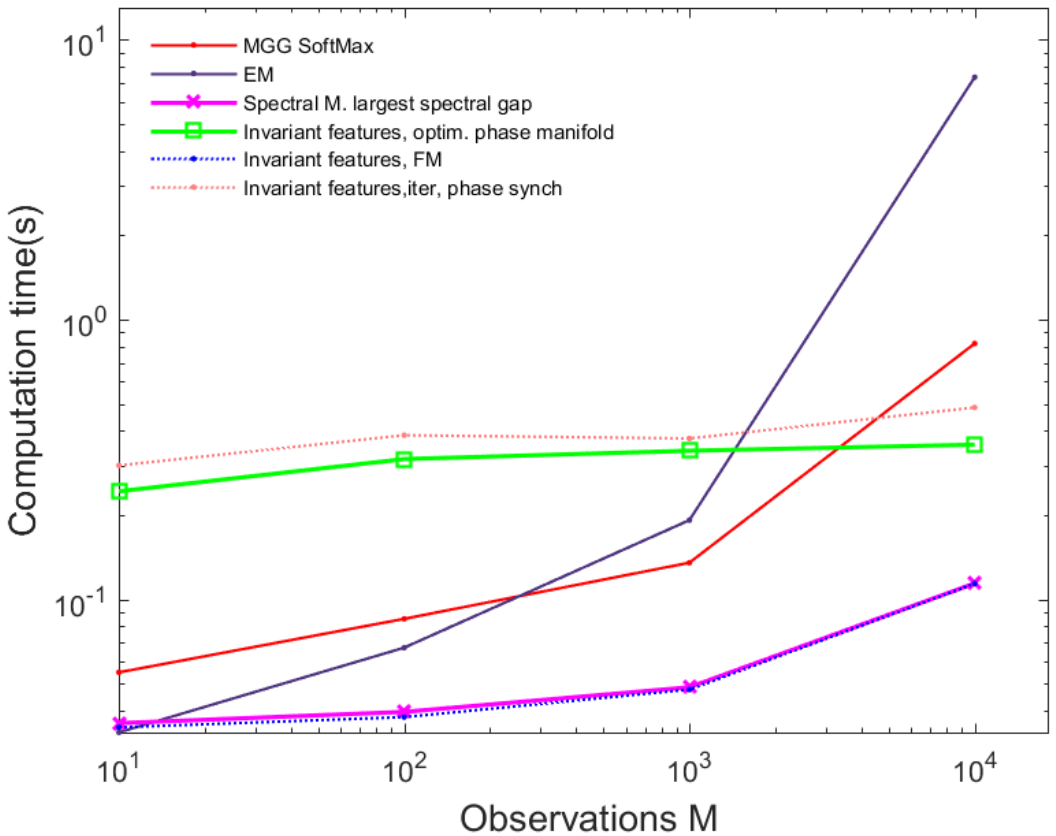}
\caption{Average computation time over 10 repetitions corresponding to \figurename\ref{Fig1}. Note that the curves corresponding to the Spectral M. largest spectral gap and FM almost overlap.}
\label{Fig2}
\end{minipage}
\end{figure}
\begin{figure}[!htb]
\begin{minipage}[t]{0.49\linewidth}
\includegraphics[width=3.15in]{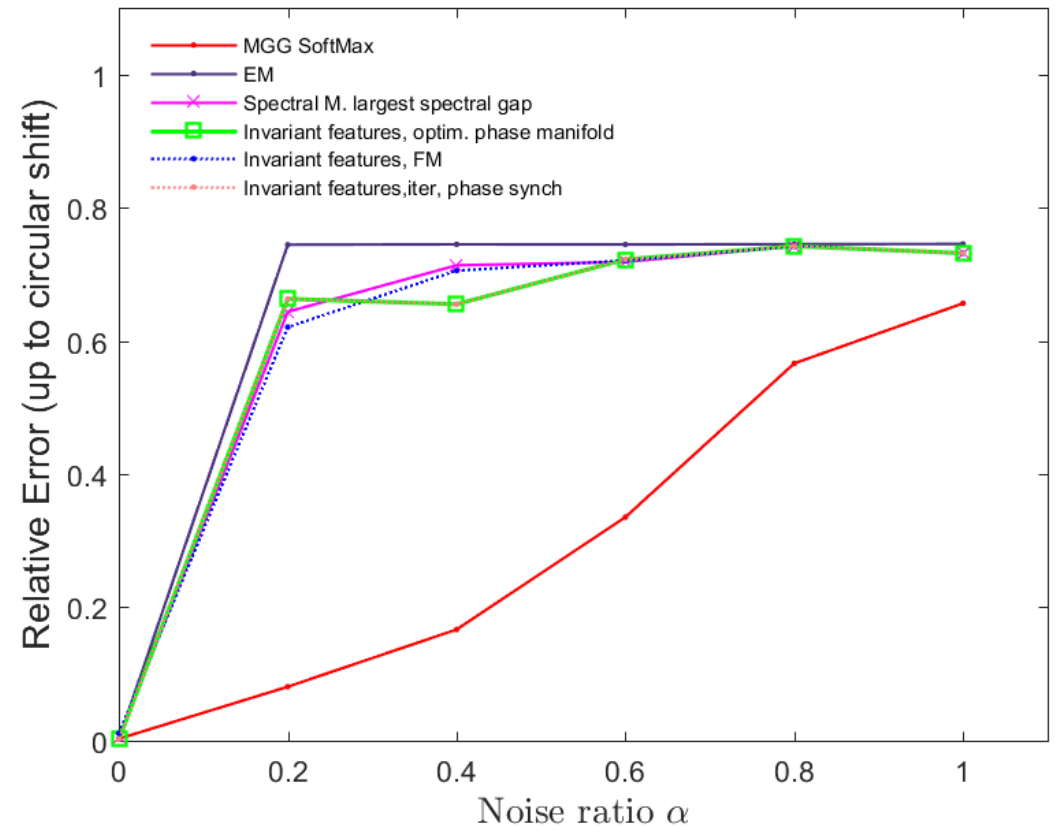}
\caption{Relative error as a function of the noise ratio $\alpha$ with fixed noise parameters $\sigma_1=10,\sigma_2=0.1$ on a real signal of length $N=41$. The number of observations is fixed at $M=10^4$. Note that curves corresponding to the optim. phase manifold and the iter. phase synch almost overlap.}
\label{Fig3}
\end{minipage}%
\vspace{0.05 in}
\hspace{0.05 in}
\begin{minipage}[t]{0.49\linewidth}
\includegraphics[width=3.15in]{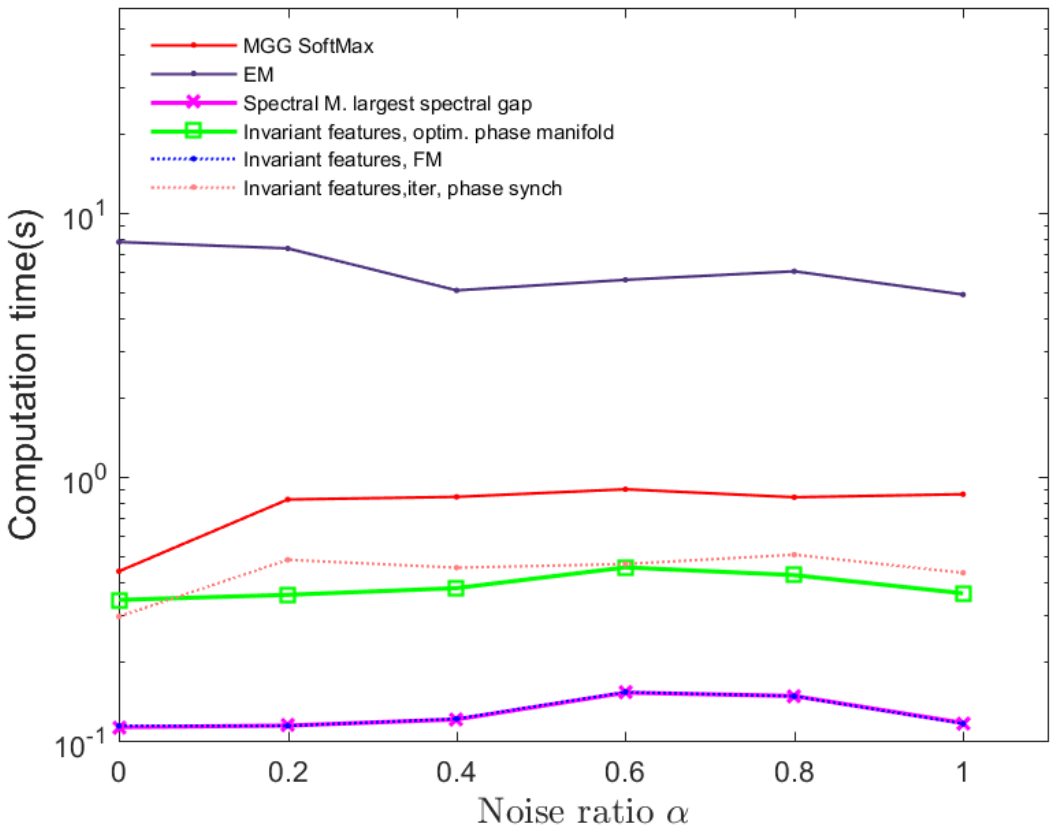}
\caption{Average computation time over 10 repetitions corresponding to \figurename\ref{Fig3}. Note that the curves corresponding to the Spectral M. largest spectral gap and FM almost overlap.}
\label{Fig4}
\end{minipage}
\end{figure}
\begin{figure}[!htb]
\begin{minipage}[t]{0.49\linewidth}
\includegraphics[width=3.15in]{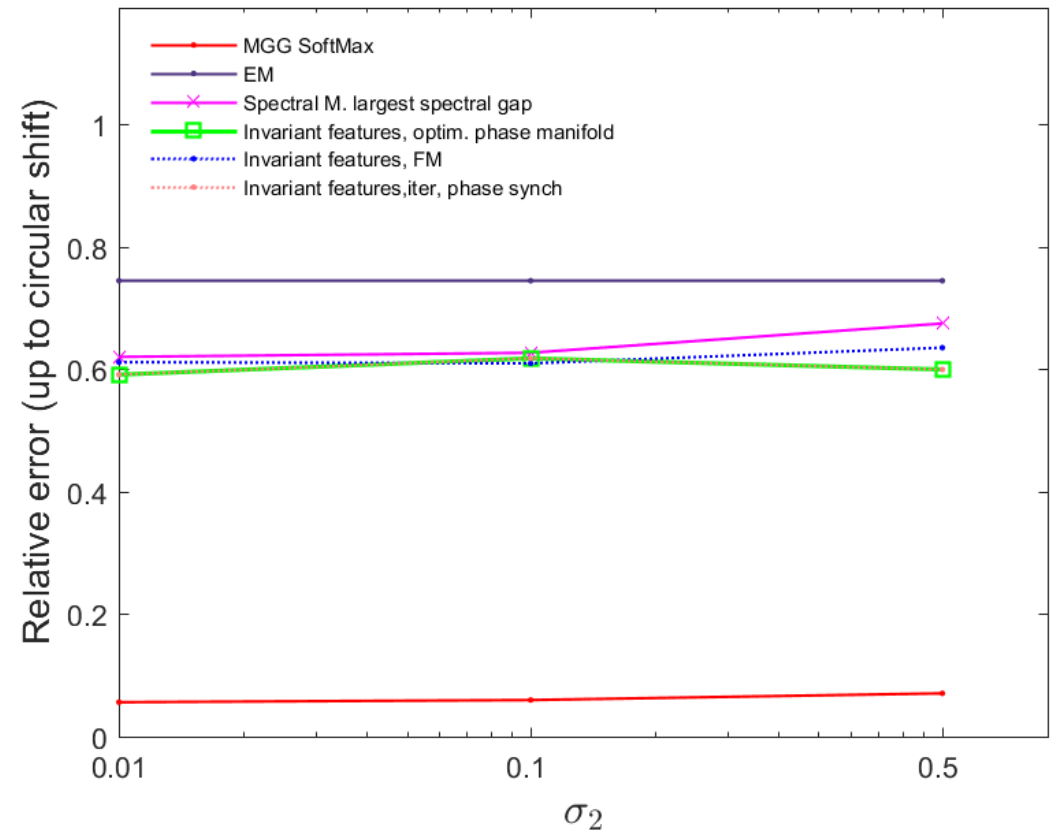}
\caption{
Relative error as a function of the second noise standard deviation $\sigma_2$ with fixed noise parameters $\alpha=0.2,\sigma_1=5$ on a real signal of length $N=41$. The number of observations is fixed at $M=10^4$. Note that curves corresponding to the optim. phase manifold and the iter. phase synch almost overlap.}
\label{Fig5}
\end{minipage}%
\vspace{0.05 in}
\hspace{0.05 in}
\begin{minipage}[t]{0.49\linewidth}
\includegraphics[width=3.15in]{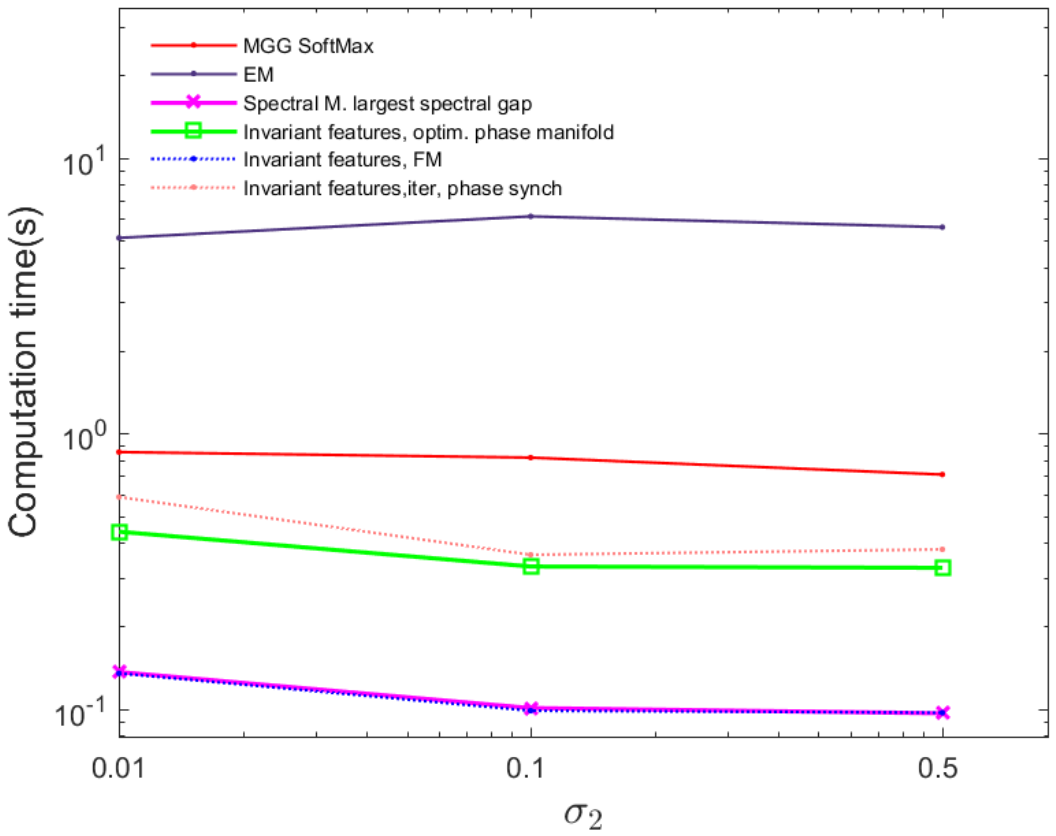}
\caption{Average computation time over 10 repetitions corresponding to \figurename\ref{Fig5}. Note that the curves corresponding to the Spectral M. largest spectral gap and FM almost overlap.}
\label{Fig6}
\end{minipage}
\end{figure}

\begin{table}[!ht]
\renewcommand\arraystretch{1.25}
  \centering
  \caption{Comparison of relative recovery error values for a fixed real piecewise constant signal
  $\pmb{u}$ of length $N=101$ under different noise levels for EM\cite{katsevich2020likelihood},
  spectral M. largest spectral gap\cite{8352518}, optim. phase manifold\cite{2017Bispectrum},
   FM\cite{2017Bispectrum}, iter. phase synch\cite{2017Bispectrum} and
    MGG SoftMax model with setting different initial values.
    (The largest relative errors are shown in bold fonts.)}
    \label{TABLE2}
  \begin{tabular}{
c@{\hspace{3mm}}c@{\hspace{3mm}}c@{\hspace{3mm}}c@{\hspace{3mm}}
c@{\hspace{3mm}}c@{\hspace{3mm}}c@{\hspace{3mm}}c@{\hspace{3mm}}c@{\hspace{3mm}}
c@{\hspace{3mm}}c@{\hspace{3mm}}c@{\hspace{3mm}}c@{\hspace{3mm}}
c@{\hspace{3mm}}c@{\hspace{3mm}}c@{\hspace{3mm}}c@{\hspace{3mm}}
c c c c c c c c c c c c c c c c c}
    \hline
\hline
\\
\hline
$\alpha$&$\sigma_1$&$\sigma_2$&&
\multicolumn{3}{c}{Existing Methods} &&&
\multicolumn{1}{c}{Proposed Method}\\
\cline{5-11}
\Xcline{4-8}{1pt}\Xcline{10-11}{1pt}
\multirow{2}*{$\downarrow$}&\multirow{2}*{$\downarrow$}&\multirow{2}*{$\downarrow$}
&\multirow{2}*{EM\cite{katsevich2020likelihood}}
&\multicolumn{4}{c}{Invariant features}&&\multirow{2}*{MGG SoftMax}
\\
\cdashline{5-8}
&&&&Spec. M.\cite{8352518}&Opt. P. manifold\cite{2017Bispectrum}&FM\cite{2017Bispectrum}
&Iter. P. synch\cite{2017Bispectrum}&&&
\\
\hline\hline
0&10&0.01
&0.0001218
&0.0157
&0.0001273
&0.0001333
&0.0001273
&&\pmb{0}
\\
0.2&10&0.01
&0.3907
&0.4563
&0.4528
&0.4523
&0.4528
&&\pmb{0}
\\
0.4&10&0.01
&0.4198
&0.6123
&0.4526
&0.4519
&0.4526
&&\pmb{0.09737}
\\
0.6&10&0.01
&0.4250
&0.6474
&0.4374
&0.4354
&0.4374
&&\pmb{0.1620}
\\
0.8&10&0.01
&0.4672
&0.5103
&0.5562
&0.4764
&0.5562
&&\pmb{0.1642}
\\
1&10&0.01
&\pmb{0.4064}
&0.4637
&0.4160
&0.4239
&0.4499
&&0.6837
\\
\hline
0&10&0.1
&0.001218
&0.1044
&0.002427
&0.002767
&0.002427
&&\pmb{0}
\\
0.2&10&0.1
&0.4520
&0.5450
&0.5452
&0.5449
&0.5452
&&\pmb{0}
\\
0.4&10&0.1
&0.5447
&0.5751
&0.5749
&0.5749
&0.5749
&&\pmb{0.1037}
\\
0.6&10&0.1
&0.4248
&0.6508
&0.4368
&0.4350
&0.4368
&&\pmb{0.1970}
\\
0.8&10&0.1
&0.4665
&0.5126
&0.5553
&0.4758
&0.5553
&&\pmb{0.2473}
\\
1&10&0.1
&\pmb{0.4064}
&0.4637
&0.4160
&0.4239
&0.4499
&&0.6837
\\
\hline
0&10&0.5
&0.006857
&0.3769
&0.02758
&0.03679
&0.02758
&&\pmb{0}
\\
0.2&10&0.5
&0.5142
&0.5676
&0.5676
&0.5674
&0.5675
&&\pmb{0.06618}
\\
0.4&10&0.5
&0.5346
&0.5705
&0.5707
&0.5705
&0.5705
&&\pmb{0.1464}
\\
0.6&10&0.5
&0.4163
&0.6563
&0.4925
&0.4252
&0.4925
&&\pmb{0.1918}
\\
0.8&10&0.5
&0.4596
&0.5223
&0.5089
&0.4680
&0.5089
&&\pmb{0.3138}
\\
1&10&0.5
&\pmb{0.4064}
&0.4637
&0.4160
&0.4239
&0.4499
&&0.6837
\\
\hline
0&5&0.01
&0.0001218
&0.0157
&0.0001273
&0.0001333
&0.0001273
&&\pmb{0}
\\
0.2&5&0.01
&0.2920
&0.4829
&0.3579
&0.3568
&0.3579
&&\pmb{0}
\\
0.4&5&0.01
&0.3712
&0.5980
&0.4685
&0.4683
&0.4686
&&\pmb{0.05013}
\\
0.6&5&0.01
&0.3213
&0.6738
&0.3539
&0.3539
&0.3539
&&\pmb{0.09212}
\\
0.8&5&0.01
&0.3371
&0.4450
&0.3583
&0.3490
&0.3583
&&\pmb{0.2541}
\\
1&5&0.01
&0.328
&0.3465
&0.3562
&0.3469
&0.3562
&&\pmb{0.3096}
\\
\hline
0&5&0.1
&0.001218
&0.1044
&0.002427
&0.002767
&0.002427
&&\pmb{0}
\\
0.2&5&0.1
&0.2925
&0.4906
&0.3573
&0.3564
&0.3574
&&\pmb{0}
\\
0.4&5&0.1
&0.3707
&0.6011
&0.4672
&0.4666
&0.4672
&&\pmb{0.06696}
\\
0.6&5&0.1
&0.3204
&0.6792
&0.3487
&0.3430
&0.3487
&&\pmb{0.1123}
\\
0.8&5&0.1
&0.3791
&0.4813
&0.4215
&0.4152
&0.4215
&&\pmb{0.2951}
\\
1&5&0.1
&0.328
&0.3465
&0.3562
&0.3469
&0.3562
&&\pmb{0.3096}
\\
\hline
0&5&0.5
&0.006857
&0.3769
&0.02758
&0.03679
&0.02758
&&\pmb{0}
\\
0.2&5&0.5
&0.2812
&0.5225
&0.3474
&0.3470
&0.3474
&&\pmb{0.06989}
\\
0.4&5&0.5
&0.3632
&0.6259
&0.4376
&0.4371
&0.4376
&&\pmb{0.1433}
\\
0.6&5&0.5
&0.3162
&0.6857
&0.3507
&0.3397
&0.3397
&&\pmb{0.1813}
\\
0.8&5&0.5
&0.3343
&0.4827
&0.3550
&0.3465
&0.3550
&&\pmb{0.2718}
\\
1&5&0.5
&0.3280
&0.3465
&0.3562
&0.3469
&0.3562
&&\pmb{0.3096}
\\
\hline
\end{tabular}
\end{table}
There are two groups of experiments with different true signals. The experiments are under the number of observations $M=10^4$. One of the true signal $\pmb{u}$ of length $N=41$ is a 1-dimensional standard Gaussian random signal. The experiments are conducted under different noise levels where the noise ratio $\alpha$ varies over $\{0,0.2,0.4,0.6,0.8,1\}$, the standard deviation of the first Gaussian noise $\sigma_1$ varies over $\{10,5\}$ and the standard deviation of the first Gaussian noise $\sigma_1$ varies over $\{0.01,0.1,0.5\}$. The corresponding relative errors are shown in \tablename~\ref{TABLE1}. The other of true signal $\pmb{u}$ of length $N=101$ is a $1$-dimensional  piecewise constant signal. Specifically, the components of $\pmb{u}$ from $30$th to  $60$th  equal to $1$, and others equal to $0$. We design the latter true signal to further show the effect of the regularization term.
The corresponding relative errors are shown in \tablename~\ref{TABLE2}. We can easily see that the relative errors in \tablename~\ref{TABLE2} are better than those in \tablename~\ref{TABLE1}. As both \tablename~\ref{TABLE1} and \tablename~\ref{TABLE2} show, it's clear to see that the proposed method named MGG SoftMax outperforms the other existing algorithms mentioned above except for the case of $\alpha=0,1$. However, notice that there is usually a tiny difference between the proposed method and the best results in the case of $\alpha=0,1$, which can be usually negligible.
\par
To intuitively see that the effect of parameters like $M,\alpha,\sigma_2$ on the relative error and computation time, we show \figurename~\ref{Fig1}-\figurename~\ref{Fig6} based on the true random signal of length $N=41$. Specific details are as follows.

\begin{description}
  \item[C1] \textbf{Effect of the number of observations ($\pmb{M}$) on the relative error and computation time.}\\
  \figurename~\ref{Fig1} and  \figurename~\ref{Fig2} show the relative error and computation time of all methods mentioned above as a function of the number of observations $M$ with fixed noise level $\alpha=0.2,\sigma_1=10,\sigma_2=0.1$. As \figurename~\ref{Fig1} shows, the proposed method named MGG SoftMax outperforms the other algorithms for $M=10,10^2,10^3,10^4$. Furthermore, with the increase of the observations, the relative error of the proposed method decreases obviously. In particular, when $M$ equals to $10^4$, the relative error of proposed method named MGG SoftMax is less than the best approaches by a factor of $65$. However, the computation time of the proposed method named MGG SoftMax is more than  the best approaches by a factor of $10$ As \figurename~\ref{Fig2} shows.
  \item[C2] \textbf{Effect of the noise ratio ($\alpha$) on the relative error and computation time.}
  \\
  We fix the number of observations $M=10^4$ and the two standard deviations $\sigma_1=10,\sigma_2=0.1$ and vary the mixture noise ratio $\alpha$ from $0$ to $1$ at intervals of $0.2$. In \figurename~\ref{Fig3}, the proposed method named MGG SoftMax outperforms than the other existing algorithms. The corresponding computation time is shown in \figurename~\ref{Fig4}.
  \item[C3] \textbf{Effect of the one of the noise standard deviation($\sigma_2$) on the relative error and computation time.}
  \\
With fixing the number of observations $M=10^4$ and noise parameters $\alpha=0.2,\sigma_1=5$,  we vary the standard deviation of the second Gaussian noise $\sigma_2$ among $\{0.01,0.1,0.5\}$. As \figurename~\ref{Fig5} shows, it's clear that the proposed method MGG softMax has a better performance with the relative error than the other existing algorithms. \figurename~\ref{Fig6} shows the corresponding computation time.
\end{description}
\section{Conclusion}\label{Sec8}
We derive a general adaptive variational model for MRA problem with Gaussian mixture noise. Compared with the existing methods, the difficulty is that besides shifts are unknown, the type of noise at each point is unknown. To overcome it, the proposed model contains two weights to determine unknown shifts and unknown parameters of the mixed noise. We prove the existence of a minimizer. Furthermore, we design an algorithm to calculate the two weights by alternate iterations separately. We provide some convergence analysis. Numerical experimental results illustrate that our model provides an impressive performance.
\par
Note that accurate estimation of two weights by updating alternately is a difficult problem in itself. In this model, the error of one step does not spread as the calculation step goes further. But the weakness is that the numerical performance of our model partly depends on initial values of mixed noise parameters. It may be useful to import other algorithms to reduce the dependence of the model on initial values.
\par
Furthermore, one can continue to apply the proposed model to some practical problems, like single particle cryo-EM problem or other scientific fields. The model of denoising and alignment of 2D images in single particle cryo-EM problem  is more complicated than that of the proposed model for 1D signal.
\section*{Acknowledgements}
This work was supported by National Natural Science Foundation of China (No.31670725) and Beijing Advanced Innovation Center for Structural Biology of Tsinghua University.

\setcounter{equation}{0}
\renewcommand\theequation{A.\arabic{equation}}
\section*{Appendix}
The proof of \pmb{Theorem} \ref{theorem1} is shown as follows.
\begin{proof}
Set
\begin{equation}\label{th1}
 \mathcal{J}_1(u,\theta,w,q)=\int_{\mathbb{X}}\int_{\Omega}w(x,y)(h_1(x,y)+h_2(x,y))dydx,
\end{equation}
where
\begin{equation}\label{th1.1}
h_1(x,y)=\int_{\Omega}\int_{\mathbb{T}}\frac{
(f(x,s)-\hat{u}(s-y))^2}{2\sigma^2(t)}
q(x,y,s,t)dtds,
\end{equation}
\begin{equation}\label{th1.2}
h_2(x,y)=\int_{\Omega}\int_{\mathbb{T}}\Big(\frac{1}{2}\log\sigma^2(t)-\log\alpha(t)\Big)
q(x,y,s,t)dtds.
\end{equation}
\begin{eqnarray}\label{th2}
 \nonumber to remove numbering (before each equation)
  \mathcal{J}_2(w,q)&=&\int_{\mathbb{X}}\int_{\Omega}\int_{\Omega}\int_{\mathbb{T}}
  w(x,y)q(x,y,s,t)\log q(x,y,s,t)dtdsdydx,\\
  \mathcal{J}_3(w)&=&\int_{\mathbb{X}}\int_{\Omega} w(x,y)\log w(x,y)dydx,\\
 \mathcal{J}_4(u) &=&\lambda J(u).
\end{eqnarray}
\par
It's easy to verify that $h_1(x,y)\geq 0,$ for any $(x,y)\in\mathbb{X}\times\Omega$. For $h_2(x,y)$, we have $h_2(x,y)\geq\int_{\mathbb{T}}|\Omega|
(\frac{1}{2}\log\sigma^2(t)-\log\alpha(t)dt)\geq |\mathbb{T}||\Omega|\left(\frac{1}{2}\log\sigma^2_{max}-\log\alpha_{min}\right)$ for any $(x,y)\in\mathbb{X}\times\Omega$. Combing the condition $0\leq w(x,y)\leq 1$ for any $ (x,y)\in\mathbb{X}\times\Omega$, $J_1(u,\theta,w,q)$ is lower bounded.
The function $z\log z$ is continuous and convex. $z\log z$ reaches its maximum value $-\frac{1}{e}$ at the point $z=e^{-1}$. So $\mathcal{J}_2(w,q)\geq -\frac{1}{e}|\mathbb{T}||\Omega|\int_{\mathbb{X}}\int_{\Omega}w(x,y)dydx\geq -\frac{1}{e}|\mathbb{T}||\Omega||\mathbb{X}|$, the last
inequality is based on $\int_{\Omega}w(x,y)dy=1$ for any $x\in\mathbb{X}$. that is to say, $\mathcal{J}_2(w,q)$ is lower bounded. At the same way, $\mathcal{J}_3(w)\geq -\frac{1}{e}|\mathbb{X}||\Omega|$. $\mathcal{J}_4(u)=\lambda J(u)dy\geq 0$. Therefore, $\mathcal{J}(u,\pmb{\theta},w,q)$ is lower bounded. Then there exists a sequence $\{(u_n,\pmb{\theta}_n,w_n,q_n)\}\subset \Gamma$ such that
\begin{equation}\label{th3}
\mathcal{J}(u_n,\pmb{\theta}_n,w_n,q_n)\rightarrow\displaystyle \inf_{(u,\pmb{\theta},w,q)\in\Gamma}\mathcal{J}(u,\pmb{\theta},w,q),~~~n\rightarrow+\infty.
\end{equation}
For any $t\in\mathbb{T}$, wen have $\pmb{\theta}_n(t)=(\sigma_n^2(t),\alpha_n(t))\subset \mathbb{K}$. Moreover, $\mathbb{K}$ is a bounded close set in $\mathbb{R}^2$. So there exists a subsequence of $\pmb{\theta}_n(t)$(still denote label as $n$) and $\pmb{\theta}(t)=(\sigma^2(t),\alpha(t))\in\mathbb{K}$ such that
\begin{equation}\label{th4}
\pmb{\theta}_n(t)=(\sigma_n^2(t),\alpha_n(t))\rightarrow\pmb{\theta}(t)=(\sigma^2(t),\alpha(t)),~~~~~n\rightarrow+\infty.
\end{equation}
That is to say, $\pmb{\theta}_n(t)\rightarrow\pmb{\theta}(t)(n\rightarrow+\infty)$ is pointwise convergence with respect of $t$.
\par
Recall that $w_n\in L^\infty(\mathbb{X}\times\Omega)$. Since $L^\infty$ is the dual space of separable linear normed space $L^1$, by Banach-alaoglu theorem, we can get that there exists a weak $*$ subsequence of $w_n$(still denote label as $n$) and weak $*$ limitation $w\in L^\infty(\mathbb{X}\times\Omega)$ such that
$w_n\rightharpoonup w(n\rightarrow+\infty)$ in $L^\infty(\mathbb{X}\times\Omega)$. That is to say, for any $\varphi\in L^1(\mathbb{X}\times\Omega)$, we have
\begin{equation}\label{th5}
\int_{\mathbb{X}}\int_{\Omega}w_n(x,y)\varphi(x,y)dydx\rightarrow \int_{\mathbb{X}}\int_{\Omega}w(x,y)\varphi(x,y)dydx, ~~~n\rightarrow+\infty.
\end{equation}
Next, we will show that $w\in\mathbb{W}$ holds, i.e., $0\leq w(x,y)\leq 1, a.e.$ in $\mathbb{X}\times\Omega$ and $\int_{\Omega}w(x,y)dy=1$ for any $x\in\mathbb{X}$. At first, we will show that $0\leq w(x,y)\leq 1~~a.e.$
in $\mathbb{X}\times\Omega$. Denote $A_1=\{(x,y)\in(\mathbb{X},\Omega)|w(x,y)>1 \},A_2=\{(x,y)\in(\mathbb{X},\Omega)|w(x,y)<0\}$. Set $\varphi_1(x,y)=\chi_{A_1}(x,y),~\varphi_2(x,y)=\chi_{A_2}(x,y)$, then $\varphi_1,\varphi_2\in L^1(\mathbb{X}\times\Omega)$. Substitute $\varphi_1$ into \eqref{th5}, then we can get
\begin{equation}\label{th6}
\iint_{A_1}w_n(x,y)dydx\rightarrow \iint_{A_1}w(x,y)dydx, ~~~n\rightarrow+\infty.
\end{equation}
If $|A_1|\neq 0$, the right side of \eqref{th6} is greater than $|A_1|$. by the sign-preserving property of limitation, there exists a large integer $N>0$ such that for any $n\geq N$, $\int_{A_1}w_n(x,y)dydx>|A_1|$ holds. However, owing to $0\leq w_n(x,y)\leq 1~a.e.~(x,y)\in\mathbb{X}\times\Omega$, the left side of \eqref{th6} equals to or less than $|A_1|$. It's a contradiction. So we can get $|A_1|=0$. At the same way, we can get $|A_2|=0$. That is to say, $0\leq w(x,y)\leq 1$ holds $a.e.$ in $\mathbb{X}\times\Omega$. Next we will show that $\int_{\Omega}w(x,y)dy=1$ holds for any $x\in\mathbb{X}$. Set $h(x)=\int_{\Omega}w(x,y)dy,~\psi(x,y)=\hat{\psi}(x)=sign(h(x)-1)$, then $\psi\in L^1(\mathbb{X}\times\Omega)$ and
\begin{equation}\label{th7}
\int_{\mathbb{X}}\int_{\Omega}w_n(x,y)\psi(x,y)dydx\rightarrow\int_{\mathbb{X}}\int_{\Omega}w(x,y)\psi(x,y)dydx,
~~~n\rightarrow+\infty,
\end{equation}
which can be rewritten as
\begin{equation}\label{th7_1}
\int_{\mathbb{X}}\hat{\psi}(x)(h_n(x)-1)dx\rightarrow\int_{\mathbb{X}}\hat{\psi}(x)(h(x)-1)dx,~~~n\rightarrow+\infty.
\end{equation}
Note that the left of \eqref{th7_1} equals to $0$ for any $n$, so we can get
\begin{equation}\label{th8}
\int_{\mathbb{X}}|g(x)-1|dx=0,
\end{equation}
which implies $g(x)=\int_{\mathbb{X}}w(x,y)dy=1,~a.e.~x\in\mathbb{X}$.
\par
Since $z\log z$ is continuous and convex, $\mathcal{J}_3(w)$ is quasiconvex and weak $*$ lower semicontinuous, i.e.
\begin{equation}\label{th9}
\lim_{n\rightarrow+\infty}\inf\int_{\mathbb{X}}\int_{\Omega}w_n(x,y)\log w_n(x,y)dydx\geq
\int_{\mathbb{X}}\int_{\Omega}w(x,y)\log w(x,y)dydx,
\end{equation}
that is to say,
\begin{equation}\label{th10}
\lim_{n\rightarrow+\infty}\inf \mathcal{J}_3(w_n)\geq \mathcal{J}_3(w).
\end{equation}
\par
By the same way, we can also get that there exists a weak $*$ subsequence of $\{q_n(x,y,s,t)\}$(still denote label as $n$) and weak $*$ limitation $q(x,y,s,t)\in\mathbb{Q}$ such that $q_n(x,y,s,t)\rightharpoonup q(x,y,s,t)~~(n\rightarrow+\infty)$, i.e., for any $\kappa\in L^1(\mathbb{T}\times\Omega)$, we have
\begin{equation}\label{th11}
\lim_{n\rightarrow+\infty}\int_{\mathbb{T}}\int_{\Omega}q_n(x,y,s,t)\kappa(s,t) dsdt\rightarrow \int_{\Omega}q(x,y,s,t)\kappa(s,t)dsdt,~~~n\rightarrow+\infty.
\end{equation}
What's more, we can get that for any $(x,y)\in\mathbb{X}\times\Omega$,
\begin{equation}\label{th12}
\lim_{n\rightarrow+\infty}\inf\int_{\mathbb{T}}\int_{\Omega}q_n(x,y,s,t)\log q_n(x,y,s,t) dsdt\geq \int_{\mathbb{T}}\int_{\Omega}q(x,y,s,t)\log q(x,y,s,t) dsdt.
\end{equation}
Denote $d_n(x,y)=\int_{\mathbb{T}}\int_{\Omega}q_n(x,y,s,t)\log q_n(x,y,s,t) dsdt,~d(x,y)=\int_{\mathbb{T}}\int_{\Omega}q(x,y,s,t)\log q(x,y,s,t) dsdt$. Then by the sign-preserving of limitation we can get for the subsequence satisfying  \eqref{th12}(still denote label as $n$), there exists a large integer $N_0$ such that for any $n\geq N_0$, $d_n(x,y)\geq d(x,y)$; moreover, Combining $0\leq w_n(x,y)\leq 1 a.e.$, for the same $n$, $w_n(x,y)(b_n(x,y)-b(x,y))\geq0$. What's more, $d_n,d\in L^1(\mathbb{X}\times\Omega)$ and
\begin{equation}\label{th13}
\begin{aligned}
&\mathcal{J}_2(w_n,q_n)-\mathcal{J}_2(w,q)=\int_{\mathbb{X}}\int_{\Omega}w_n(x,y)d_n(x,y)dydx-
\int_{\mathbb{X}}\int_{\Omega}w(x,y)d(x,y)dydx\\
=&\int_{\mathbb{X}}\int_{\Omega}w_n(x,y)(d_n(x,y)-d(x,y))dydx+
\int_{\mathbb{X}}\int_{\Omega}(w_n(x,y)-w(x,y))d(x,y)dydx.
\end{aligned}
\end{equation}
By Fatou's lemma, we can get
\begin{equation}\label{th14}
\begin{aligned}
&\lim_{n\rightarrow+\infty}\inf\int_{\mathbb{X}}\int_{\Omega} w_n(x,y)(d_n(x,y)-d(x,y))dydx\\
\geq&
\int_{\mathbb{X}}\int_{\Omega} \lim_{n\rightarrow+\infty}\inf w_n(x,y)(d_n(x,y)-d(x,y))dydx\\
\geq&\int_{\mathbb{X}}\int_{\Omega} \lim_{n\rightarrow+\infty}\inf w_n(x,y)\lim_{n\rightarrow+\infty}\inf(d_n(x,y)-d(x,y))dydx\\
\geq&0.\\
\end{aligned}
\end{equation}
Substitute $\varphi(x,y)=d(x,y)$ into \eqref{th5}, we can get
\begin{equation}\label{th15}
\lim_{n\rightarrow+\infty}\int_{\mathbb{X}}\int_{\Omega}(w_n(x,y)-w(x,y))d(x,y)dydx
=0.
\end{equation}
Combing \eqref{th13},\eqref{th14} and \eqref{th15}, then
\begin{equation}\label{th16}
\lim_{n\rightarrow+\infty}\inf \mathcal{J}_2(w_n,q_n)\geq \mathcal{J}_2(w,q).
\end{equation}
\par
Since $\displaystyle\lim_{u\rightarrow+\infty}\mathcal{J}(u,\pmb{\theta},w,q)=+\infty$ and the sequence $\{(u_n,\pmb{\theta}_n,w_n,q_n)\}\subset\Gamma$ is a minimizing sequence of $\mathcal{J}(u,\pmb{\theta},w,q)$, we can get
the sequence $\{u_n\}$ is uniformly bounded with respect to $n$ and $x$ and the boundary is denoted as $M_u$. By the definition of the sequence $\{(u_n,\theta_n,w_n,q_n)\}$, then when $n$ is large enough, $\mathcal{J}(u_n,\theta_n,w_n,q_n)$ is bounded, which denoted by $C$,i.e. $\mathcal{J}(u_n,\theta_n,w_n,q_n)\leq C$. Recall that $\mathcal{J}_1,\mathcal{J}_2,\mathcal{J}_3$ is lower bounded, then $\mathcal{J}_4(u_n)$ is upper bounded. Recall $\|u_n\|_{BV(\Omega)}=\|u_n\|_{L^1(\Omega)}+J(u_n)$, then $\{u_n\}$ is bounded in $BV(\Omega)$.
So there exists a subsequence(label still denoted as $n$) and $u\in BV(\Omega)$ such that $u_n\rightarrow u$ strongly in $BV(\Omega)$ and $Du_n\rightharpoonup Du$ in the sense of distribution,i.e.,$\langle Du_n,\psi\rangle\rightarrow\langle Du,\psi\rangle,~(n\rightarrow+\infty)$ for all $\psi\in (C_0^\infty(\Omega))^2$. By the lower semicontinuity of the total variation and Fatou's lemma, we can get
\begin{equation}\label{th17}
\lim_{n\rightarrow+\infty}\inf \mathcal{J}_4(u_n)\geq \mathcal{J}_4(u).
\end{equation}
 Due to $u_n(y)\geq 0~a.e.$ and $u_n\rightarrow u$ in $L^1(\Omega)$ strongly, we can get $u\geq 0~a.e.$. So we can get $u\in S(\Omega)$.
\par
\par
Now let's consider the first term $\mathcal{J}_1(u,\pmb{\theta},w,q)$. Recall that
\begin{equation*}
h_1(x,y)=\int_{\Omega}\int_{\mathbb{T}}\frac{
(f(x,s)-\hat{u}(s-y))^2}{2\sigma^2(t)}
q(x,y,s,t)dtds,
\end{equation*}
then
\begin{equation}\label{th18}
\begin{array}{rl}
&h_{1,n}(x,y)-h_1(x,y)\\
=&\displaystyle\int_{\Omega}\int_{\mathbb{T}}\frac{
(f(x,s)-\hat{u}_n(s-y))^2}{2\sigma_n^2(t)}
q_n(x,y,s,t)dtds-\int_{\Omega}\int_{\mathbb{T}}\frac{
(f(x,s)-\hat{u}(s-y))^2}{2\sigma^2(t)}
q(x,y,s,t)dtds\\
=&\displaystyle\int_{\Omega}\int_{\mathbb{T}}\left(\frac{
(f(x,s)-\hat{u}_n(s-y))^2}{2\sigma_n^2(t)}-\frac{
(f(x,s)-\hat{u}(s-y))^2}{2\sigma^2(t)}\right)q_n(x,y,s,t)dtds\\
&\displaystyle+\int_{\Omega}\int_{\mathbb{T}}\frac{
(f(x,s)-\hat{u}(s-y))^2}{2\sigma^2(t)}
(q_n(x,y,s,t)-q(x,y,s,t))dtds\\
\end{array}
\end{equation}
For any fixed point $(s,t)\in\mathbb{T}\times\Omega$, substitute $\kappa(s,t)=\frac{
(f(x,s)-\hat{u}(s-y))^2}{2\sigma^2(t)}\in L^1(\mathbb{T}\times\Omega)$ into \eqref{th11}, then we can get
\begin{equation}\label{th19}
\lim_{n\rightarrow+\infty}\int_{\Omega}\int_{\mathbb{T}}\frac{
(f(x,s)-\hat{u}(s-y))^2}{2\sigma^2(t)}
(q_n(x,y,s,t)-q(x,y,s,t))dtds=0.
\end{equation}
\begin{equation}\label{th20}
\begin{array}{rl}
&\displaystyle\int_{\Omega}\int_{\mathbb{T}}\left(\frac{
(f(x,s)-\hat{u}_n(s-y))^2}{2\sigma_n^2(t)}-\frac{
(f(x,s)-\hat{u}(s-y))^2}{2\sigma^2(t)}\right)q_n(x,y,s,t)dtds\\
=&\displaystyle\int_{\Omega}\int_{\mathbb{T}}\left(\frac{
(f(x,s)-\hat{u}_n(s-y))^2}{2\sigma_n^2(t)}-\frac{
(f(x,s)-\hat{u}(s-y))^2}{2\sigma_n^2(t)}\right)q_n(x,y,s,t)dtds\\
&\displaystyle+\int_{\Omega}\int_{\mathbb{T}}\left(\frac{
(f(x,s)-\hat{u}(s-y))^2}{2\sigma_n^2(t)}-\frac{
(f(x,s)-\hat{u}(s-y))^2}{2\sigma^2(t)}\right)q_n(x,y,s,t)dtds
\end{array}
\end{equation}
Since $\left|(f(x,s)-\hat{u}(s-y))^2\left(\frac{1}{2\sigma_n^2(t)}
-\frac{1}{2\sigma^2(t)}\right)q_n(x,y,s,t)\right|\leq\frac{2(f(x,s)-\hat{u}(s-y))^2}{\sigma_{min}^2}\in L^1(\mathbb{T}\times\Omega),$ by Lebesgue Control Convergent Theorem, we can get
\begin{equation}\label{th21}
\begin{array}{rl}
&\displaystyle\lim_{n\rightarrow+\infty}\int_{\Omega}\int_{\mathbb{T}}(f(x,s)
-\hat{u}(s-y))^2\left(\frac{1}{2\sigma_n^2(t)}
-\frac{1}{2\sigma^2(t)}\right)q_n(x,y,s,t)dtds\\
=&\displaystyle\int_{\Omega}\int_{\mathbb{T}}(f(x,s)
-\hat{u}(s-y))^2\lim_{n\rightarrow+\infty}\left(\frac{1}{2\sigma_n^2(t)}
-\frac{1}{2\sigma^2(t)}\right)q_n(x,y,s,t)dtds\\
=&0,
\end{array}
\end{equation}
where the last equation is based on the facts $\displaystyle\lim_{n\rightarrow+\infty}\left(\frac{1}{2\sigma_n^2(t)}
-\frac{1}{2\sigma^2(t)}\right)=0$ and $0\leq q_n(x,y,s,t)\leq 1$ is bounded in $\Gamma$ for any $n$.
\begin{equation}\label{th22}
\begin{array}{rl}
&\displaystyle\left|\lim_{n\rightarrow+\infty}\int_{\Omega}\int_{\mathbb{T}}\left(\frac{
(f(x,s)-\hat{u}_n(s-y))^2}{2\sigma_n^2(t)}-\frac{
(f(x,s)-\hat{u}(s-y))^2}{2\sigma_n^2(t)}\right)q_n(x,y,s,t)dtds\right|\\
\leq&\displaystyle\lim_{n\rightarrow+\infty}\int_{\Omega}\int_{\mathbb{T}}\left(\frac{
|2f(x,s)-\hat{u}_n(s-y)-\hat{u}(s-y)||\hat{u}(s-y)-\hat{u}_n(s-y)|}{2\sigma_n^2(t)}\right)q_n(x,y,s,t)dtds\\
\leq&\displaystyle \frac{\sup f+M_u}{\sigma_{min}^2}\lim_{n\rightarrow+\infty}\int_{\Omega}\int_{\mathbb{T}}|\hat{u}(s-y)-\hat{u}_n(s-y)|dsdt
\\
=&0.
\end{array}
\end{equation}
where the last equation is based on $u_n\rightarrow u,~(n\rightarrow+\infty)$ strongly in $L^1(\mathbb{T}\times\Omega)$.
Combing \eqref{th18}, \eqref{th19}, \eqref{th20}, \eqref{th21} and \eqref{th22}, then for any $(x,y)\in\mathbb{X}\times\Omega$,
\begin{equation}\label{th23}
\lim_{n\rightarrow+\infty}h_{1,n}(x,y)=h_1(x,y).
\end{equation}
By the similar way, we can easily get
\begin{equation}\label{th23_1}
\lim_{n\rightarrow+\infty}h_{2,n}(x,y)=h_2(x,y).
\end{equation}
Furthermore, for any fixed $(x,y)\in\mathbb{X}\times\Omega$, we have
\begin{equation}\label{th24}
\lim_{n\rightarrow+\infty}h_{n}(x,y)=\lim_{n\rightarrow+\infty}h_{1,n}(x,y)
+\lim_{n\rightarrow+\infty}h_{2,n}(x,y)
=h(x,y).
\end{equation}
Next we consider
\begin{equation}\label{th25}
\begin{array}{rl}
&\mathcal{J}_1(u_n,\theta_n,w_n,q_n)-\mathcal{J}_1(u,\theta,w,q)\\
=&\int_{\mathbb{X}}\int_{\Omega}(w_n(x,y)h_n(x,y)-w(x,y)h(x,y))dydx\\
=&\int_{\mathbb{X}}\int_{\Omega}w_n(x,y)(h_n(x,y)-h(x,y))dydx
+\int_{\mathbb{X}}\int_{\Omega}(w_n(x,y)-w(x,y))h(x,y)dydx.
\end{array}
\end{equation}
Substitute $\varphi(x,y)=h(x,y)\in L^1(\mathbb{X}\times\Omega)$ into \eqref{th5}, we can get
\begin{equation}\label{th26}
\int_{\mathbb{X}}\int_{\Omega}(w_n(x,y)-w(x,y))h(x,y)dydx=0.
\end{equation}
 By Lebesgue Control Convergence Theorem, we can get
\begin{equation}\label{th27}
\begin{aligned}
&\displaystyle\left|\lim_{n\rightarrow+\infty}\int_{\mathbb{X}}\int_{\Omega}w_n(x,y)(h_n(x,y)-h(x,y))dydx\right|\\
\leq&\displaystyle\left|\lim_{n\rightarrow+\infty}\|w_n\|_{L^\infty(\mathbb{X}\times\Omega)}
\int_{\mathbb{X}}\int_{\Omega}(h_n(x,y)-h(x,y))dydx\right|\\
\leq&\displaystyle\left|\lim_{n\rightarrow+\infty}
\int_{\mathbb{X}}\int_{\Omega}(h_n(x,y)-h(x,y))dydx\right|\\
=&0.
\end{aligned}
\end{equation}
The last equality sign is based on \eqref{th24}. So
\begin{equation}\label{th28}
\lim_{n\rightarrow+\infty}\int_{\mathbb{X}}\int_{\Omega}w_n(x,y)(h_n(x,y)-h(x,y))dydx=0.
\end{equation}
Substitute \eqref{th26} and \eqref{th28} into \eqref{th25}, we can get
\begin{equation}\label{th29}
 \lim_{n\rightarrow+\infty}\mathcal{J}_1(u_n,\pmb{\theta}_n,w_n,q_n)=\mathcal{J}_1(u,\pmb{\theta},w,q).
\end{equation}
Above all, we can get
\begin{equation}\label{th30}
\lim_{n\rightarrow+\infty}\mathcal{J}(u_n,\pmb{\theta}_n,w_n,q_n)\geq \mathcal{J}(u,\pmb{\theta},w,q),
\end{equation}
which implies $(u,\pmb{\theta},w,q)$ is a minimizer.
\end{proof}
\end{document}